\documentclass[12pt,a4paper,twoside,final,notitlepage, reqno]{article}

\usepackage[latin1]{inputenc}
\usepackage[english]{babel}
\usepackage[a4paper,left=2.5cm,right=2.5cm]{geometry}
\usepackage{amssymb}
\usepackage{amsthm} 
\usepackage{amsmath}
\usepackage{graphics,graphicx}
\usepackage{epstopdf}
\usepackage[usenames, dvipsnames]{color}
\usepackage[textsize=small]{todonotes}
\usepackage{booktabs}
\usepackage{subfigure} 
\usepackage{siunitx}

\graphicspath{{figs/}}

\makeatletter
\def\input@path{{figs/}}
\makeatother

\theoremstyle{definition}
\newtheorem{theorem}{Theorem}%[section]
\newtheorem{lemma}{Lemma}
\newtheorem{proposition}{Proposition}

\theoremstyle{definition}
\theoremstyle{definition}
\newtheorem{remark}{Remark}%[section]

\providecommand{\keywords}[1]  {\textbf{Keywords:} #1}
\providecommand{\subjclass}[1] {\textbf{Subject classification:} #1}
\providecommand{\acknow}[1] {\textbf{Acknowledgments.} #1}

\def\sizeg{ \footnotesize}

\usepackage{esdiff}
\renewcommand{\d}[1]{\, \mathrm{d} #1}

\DeclareMathOperator{\diag}{diag}

\DeclareMathOperator{\e}{e}

\newcommand{\C}{\mathbb{C}}
\newcommand{\R}{\mathbb{R}}

\renewcommand{\phi}{\varphi}
\renewcommand{\epsilon}{\varepsilon}

\def \th{_{t,h}}

\def \RL{  \text{RL} }
\def \EAB{ \text{EAB} }

\usepackage{authblk}

\title{
  \vspace{-2.5cm}
  \bf{\Large{
      Rush-Larsen time-stepping methods of high order 
      \\[5pt] 
      for stiff problems in cardiac electrophysiology
    }}
}

\author[1,2]{Yves Coudi\`ere \thanks{yves.coudiere@inria.fr}}
\author[1]{Charlie Douanla-Lontsi \thanks{charlie.douanla-lontsi@inria.fr}}
\author[3]{Charles Pierre \thanks{charles.pierre@univ-pau.fr}}

\affil[1]{
  INRIA Bordeaux Sud Ouest, Universit\'e de Bordeaux, France.
}
\affil[2]{
  Institut de Math\'ematiques de Bordeaux, 
  UMR CNRS 5241.
}\affil[3]{
  Laboratoire  de Math\'ematiques et de leurs Applications,
  UMR CNRS 5142, \protect \\
  Universit\'e de Pau et des Pays de l'Adour, France.
}

\usepackage{fancyhdr}

\begin{document} 

\date{15 October, 2018}

\maketitle

\begin{abstract}
  Stability and accuracy of numerical methods for reaction-diffusion equations
  still need improvements, which prompts for the development of high order and
  stable time-stepping methods.
  This is particularly true in the context of cardiac
  electrophysiology, where reaction-diffusion equations are coupled with stiff
  systems of ordinary differential equations. 
  So as to address
  these issues, much research on implicit-explicit methods and exponential
  integrators has been carried out during the past 15 years. In 2009,
  Perego and Veneziani~\cite{perego-2009} proposed an innovative
  time-stepping scheme of order 2. In this paper we present an
  extension of this scheme to the orders 3 and 4, that we call
  Rush-Larsen schemes of order $k$. These new schemes are
  explicit multistep methods, which belong to the classical class of
  exponential integrators. Their general formulation is simple and easy
  to implement. We prove that they are stable under perturbation and
  convergent of order $k$. We analyze their Dahlquist stability,
  and show that they have a very large stability domain, provided that
  the stabilizer associated with the method captures well enough the stiff modes
  of the problem. We study their application to a system of equations
  that models the action potential in cardiac electrophysiology.
\end{abstract}
\vspace{20pt}
\noindent
\\[3pt]
\\[3pt]
\acknow{
  This study received financial support from the French government as part
  of the ``Investissement d'Avenir'' program managed by the ``Agence Nationale de
  la Recherche'' (ANR), grant reference ANR-10-IAHU-04. It also received fundings
  of the ANR project HR-CEM, grant reference 13-MONU-0004-04.
}
\\[20pt]
\section*{Introduction}
This article concerns the problem of time integration of stiff
reaction-diffusion equations, in particular 
when they are coupled to a system
of ordinary differential equations (ODE). As developed below, for such problems, the matters of stability and accuracy are of
first importance. As a systemic example of these questions, we will consider the
monodomain model in cardiac electrophysiology \cite{CLEMENS-NENONEN-04,
  FRANZONE-PAVARINO-TACCARDI-05.1,
  FRANZONE-PAVARINO-TACCARDI-05.2}. Given the heart domain $\Omega$ and
the time interval $[0,T]$, it has the general form
\begin{equation}
  \label{F0}
  \dfrac{\partial v}{\partial t} = Av + f_1(v,\zeta) + s(x,t), \quad 
  \diff{\zeta}{t} = f_2(v,\zeta),
\end{equation}
where $A$ is a diffusion operator. The unknown function
$v:\Omega\times [0,T]\to \R$ is the transmembrane potential. 
The unknown
function $\zeta:\Omega\times [0,T]\to\R^{p+q}$ gathers $p+q$ variables
describing the state of the cell membrane. It incorporates $p$ gating variables
and $q$ ionic concentrations. The source term $s(x,t)$ is an applied stimulation
current. The reaction terms $f_1$ and $f_2$ model ionic currents across the cell
membrane, and are called ionic models. Ionic models have originally been
developed by Hodgkin and Huxley \cite{Hodgkin52} in 1952. Highly detailed
ionic models specific to cardiac cells have been designed since the 1960's,
such as the Beeler and Reuter (BR) model \cite{beeler-reuter}
% , the Luo
% and Rudy \blue{(LR)} models \cite{LR1,LR2a}
or the ten Tusscher, Noble,
Noble and Panfilov (TNNP) model \cite{tnnp}. A comprehensive review is
available in \cite{spiteri}.

There are two major difficulties for numerical simulations in cardiac
electrophysiology. First, the non-linear functions $f_1$ and $f_2$ in
equation~\eqref{F0} induce expensive computations of the mappings
$(v,\zeta)\rightarrow f_i(v,\zeta)$. For example, the TNNP model \cite{tnnp}
involves the computation of 50 scalar exponentials, that have to
be performed for each mesh node to approximate solutions of the partial
differential equation~\eqref{F0}. They represent the predominant
computational load during numerical 
simulations, and their
total amount needs to be
maintained as low as possible. Fully implicit time-stepping methods, which
require a non-linear solver, are therefore avoided.
Second, the equations~\eqref{F0} are stiff, 
but since implicit methods are not
affordable, numerical instabilities are
challenging to manage.
More precisely, the stiffness is caused by the presence of
different space and time scales. The solutions of equation~\eqref{F0}
display sharp wavefronts. Typically, the scaling factor between the fast
and slow variables ranges from 100 to 1000. This is commonly coped with by
resorting to very fine space and time discretization grids, associated with high
computational costs.

In this context, our strategy for solving problem~\eqref{F0} is to use high
order methods, so as to have accurate simulations with coarser space and
time discretization grids. A high order time-stepping method that fulfills the
two following conditions is required: it must have strong stability
properties, and has to be explicit for the reaction terms. To this aim, we will focus
on the time integration of stiff ODE systems of the form
\begin{equation}
  \label{F1}
  \diff{y}{t} = f(t,y), \quad y(0) = y_{0},
\end{equation}
in which the nonlinear function $f: [0,+\infty[\times \R^N \rightarrow \R^N$ (e.g. $N=p+q$
for the ionic models presented above) may be split as
$f(t,y) = a(t,y)y+b(t,y)$. This leads to a formulation more suited to
our needs,
\begin{equation}
  \label{F2}
  \diff{y}{t} = a(t,y)y+b(t,y), \quad y(0) = y_{0}.
\end{equation} 
It involves the non-linear term $b(t,y)$ and the operator
$y \in\R^N \mapsto a(t,y)y \in\R^N$, which can be easily linearized as
e.g. $a(t,y_n)y$.  This term $a(t,y)$ will be inserted into the
numerical scheme in order to stabilize the computations. It will be called
\textit{the stabilizer} in the sequel. In practice, the term
$a(t,y)$ may be related to the Jacobian of the system
$\partial _y f(t,y)$. However, no a priori definition of the stabilizer
is made (such as $a(t,y)=\partial _y f(t,y)$), because we plan to
analyze the formulation in~\eqref{F2} in general. 
This will allow us  for instance
to define the stabilizer as an approximation
of the Jacobian,
for technical reasons detailed below.
This approach is relevant in cardiac electrophysiology, where
the fastest variables are gating variables that
are given by the $p$ first equations
of the ODE system $\d{\zeta}/\d{  t} = f_{2}(v,\zeta)$
in (\ref{F0}). 
They have the
general form 
\begin{equation*}
  \diff{  \zeta_i}{  t} =  
  % a_i(v)w_i + b_i(v),
  \dfrac{\zeta_{i,\infty} (v) - \zeta_i}{\tau_i(v)},
\end{equation*}
(see Section \ref{sec:memb-eq}) that
motivates the reformulation \eqref{F2} with the diagonal stabilizer $a=\diag(-1/\tau_i)$.

Exponential integrators are well suited in this 
framework, we  refer to
\cite{Bor_W, Hochbruck-Ostermann-2010, hochbruck-2015} for general reviews. They
have been widely studied for the semilinear equation
$\partial _t y = Ay + b(t,y) $, see \textit{e.g.}~\cite{hochbruck-98,
  Cox_Mat,oster_ex_RK, exp-rozenbrock-2009, Tokman-2012, Ostermann-ERK-2014}.
Exponential integrators commonly define a time iteration based on the
exact solution of an equation of the form $ \partial _t y = Ay + p(t) $ where
$p(t)$ is a polynomial. It is usually defined with the functions
$\left(\phi_k\right)_{k\ge 0}$
\begin{equation}
  \label{eq:phi-j}
  \phi_0(z) = \e^z,\quad \phi_{j+1}(z) = \frac{\phi_{j}(z) - 1/j!}{z},
\end{equation}
introduced by N{\o}rsett \cite{Nors_EAB}.
% The basic idea is to use the exact solution of the linearized equation in
% order to stabilize the numerical scheme.
In general, it requires to compute a matrix exponential applied to a
vector, like $\e^{tA}y$.%  \blue{The difficulty \blue{associated with} the cost of
% this operation is largely mitigated when the matrix is diagonal.}  \question{à
% mettre ou pas à la place de l'autre phrase ?} 
This is the supplementary cost
associated with exponential integrators. A gain in stability is expected
when $A$ is the predominant stiff part of the equation.
% , specifically when
% $b(t,y)$ has locally a smaller Lipschitz constant than $Ay + b(t,y)$.}
% \question{Je sais pas s'il faut laisser la dernière partie, derrière la virgule,
% qui n'est pas très claire pour moi.}

The target equation~\eqref{F2} incorporates a non-constant
linear part $a(t,y)$, exponential integrators have been less studied in that
case.  Exponential integrators of Adams type for a non-constant linear part have
been first considered by Lee and Preiser \cite{Lee_Stan} in 1978, and by Chu
\cite{CHU} in 1983. Recently, Ostermann et al., \cite{ EAB_M_O} developed and
analyzed the linearized exponential Adams method. In general, the original
equation~\eqref{F1} is formulated after each time step as
$\frac{\d{y}}{\d{t}} = J_n y + c_n(t,y)$,
% \begin{equation*}
%   \diff{y}{t} = J_n y + c_n(t,y), %\quad J_n = \partial _y f(t_n,y_n), \quad 
%   %   c_n(t,y) = f(t,y) - J_n y,
% \end{equation*} 
involving the Jacobian matrix $J_n= \partial _y f(t_n,y_n)$, and the
correction function $c_n(t,y) = f(t,y) - J_n y$. This has several drawbacks.
It requires the computation of matrix exponential applied to a vector with a
different matrix at each time step. Moreover, 
stabilization can be performed on the fast variables only, in case they are
known in advance, e.g. because of modeling assumptions, or of our physical
understanding of the problem. In this case, using the full Jacobian as the
stabilizer will cause unnecessary computational efforts. As an alternative, the
stabilizer can be set to a part or an approximation of the Jacobian.
This had already been proposed by N{\o}rsett \cite{Nors_EAB} in 1969
and has been analyzed in \cite{Tranquilli-Sandu-2014},
\cite{Rainwater-Tokman-2016}, and \cite{coudiere-lontsi-pierre-2015} for
exponential Rosenbrock, exponential Runge-Kutta and exponential Adams methods,
respectively. For exponential Adams methods, equation~\eqref{F2} is reformulated
after each time step as $\frac{\d{y}}{\d{t}} = a_n y + c_n(t,y)$, with
$a_n = a(t_n,y_n)$, and $c_n(t,y) = f(t,y) - a_n y$.
% \begin{displaymath}
%   \diff{y}{t} = a_n y + c_n(t,y), \quad a_n = a(t_n,y_n), \quad 
%   c_n(t,y) = f(t,y) - a_n y.
% \end{displaymath}  
The resulting scheme with a time-step $h>0$ is (see the details in
\cite{EAB_M_O,coudiere-lontsi-pierre-2015})
\begin{equation}
  \label{eab-scheme}
  y_{n+1} = y_n + h \left( \phi_1(a_n h) \left( a_n y_n + \gamma_{1} \right) +
    \phi_2(a_n h)\gamma_{2} + \ldots \phi_k(a_n h)\gamma_{k} \right),
\end{equation}
where the numbers $\gamma_{i}$ are the coefficients of the Lagrange
interpolation polynomial of $c_n(t,y)$ (in a classical $k$-step setting), and
the functions $\phi_j$ are given by~\eqref{eq:phi-j}.

Independently, Perego and Veneziani \cite{perego-2009} presented in 2009 an
innovative exponential integrator of order 2, of a different nature. They
proposed a scheme of the form
\begin{equation}
  \label{eq:RL}
  y_{n+1} = y_n + h \phi_1(\alpha_n h) \left( \alpha_n y_n + \beta_n \right),
\end{equation}
involving two coefficients $\alpha_n$ and $\beta_n$ to be computed at
each time step.  The resulting scheme has a very simple 
definition, and is in
particular simpler than the exponential Adams integrators~\eqref{eab-scheme}.
The essential difference with the previous approaches is that
$\alpha_n \ne a(t_n,y_n)$, but instead is fixed for the scheme to be
consistent of order 2. Specifically, the coefficients
$\alpha_n$ and $\beta_n$ are given by $\alpha_n = \frac32 a_n - \frac12 a_{n-1}$
and $ \beta_n = \frac32 b_n - \frac12 b_{n-1}$ with $a_j = a(t_j, y_j)$ and
$b_j = b(t_j, y_j)$. Perego and Veneziani presented their scheme as a
\textit{``generalization of the Rush-Larsen method''} in reference to the
Rush-Larsen scheme \cite{RL1} commonly used in electrophysiology. 
% \blue{Indeed,
% the Rush-Larsen scheme has the form~\eqref{eq:RL} with $\alpha_n = a_n$, which
% guarantees a consistency of order 1.} \question{correct ?}

This scheme resembles the Magnus integrator introduced by Hochbruck et
al. in \cite{hochbruck2003magnus} for the time dependent Schrödinger equation
$i y' = H(t)y$, and extended by Gonzàlez et al. in \cite{gonzalez2006second}
to parabolic equations with time-dependent linear part $y' = a(t)y + b(t)$.
The second-order Magnus integrator also formulates as~\eqref{eq:RL}, but
with $\alpha_n=a(t_{n+1/2})$ and $b_n=b(t_{n+1/2})$. The scheme of Perego and
Veneziani generalizes the second-order Magnus integrator to the case where
$a = a(t,y)$ and $b = b(t,y)$: it presents an approximation of the unknown terms
$a\left( t_{n+1/2}, y(t_{n+1/2}) \right)$ and
$b\left( t_{n+1/2}, y(t_{n+1/2}) \right)$ using a two-points interpolation.

In this paper we will study schemes under the form~\eqref{eq:RL}.  We
will show that they also exist at the orders 3 and 4, and will
exhibit explicit definitions of the two coefficients $\alpha_n$ and
$\beta_n$. The schemes will be referred to as as Rush-Larsen schemes of
order $k$ (shortly denoted by $\RL_k$), in the continuation of the
denomination used in \cite{perego-2009}. They will be shown to be stable under
perturbation and convergent of order $k$. We also present the Dahlquist
stability analysis for the $\RL_k$ schemes.  It is a practical tool that
allows one to scale the time step $h$ with respect to the
variations of the function $f(t,y)$ in problem~\eqref{F1}, see
\textit{e.g.}  \cite{Hairer-ODE-II}. The splitting
$f(t,y) = a(t,y)y+b(t,y)$ may be arbitrary, but obviously the choice of the
stabilizer term $a(t,y)$ is critical for the stability of the method. When
considering time-dependent stabilizers, the stability domain depends on
this splitting. We compute stability domains numerically, and show that they are
much larger if $a(t,y)$ captures the variations of $f(t,y)$, than in absence of
stabilization (i.e., when $a(t,y)=0$). We finally evaluate the
performances of the $\RL_k$ methods as applied to the membrane equation in cardiac electrophysiology. They are compared to the exponential
Adams integrators~\eqref{eab-scheme}. The two methods have a very similar
robustness with respect to stiffness, allowing stable
computations with large time steps. 
For the considered test case,
the
$\RL_3$ and $\RL_4$ schemes are more accurate for large time steps.

The paper is organized as follows. The $\RL_k$ schemes are derived in Section
\ref{sec:RL-def}, and their numerical analysis is made in Sections
\ref{sec:RL-def} and \ref{sec:conv}. The Dahlquist stability analysis is
completed in Section \ref{sec:dahlquist-stab}. The numerical results are
presented in Section \ref{sec:num-res}. The paper ends with a conclusion
in Section \ref{sec:conclusion}.

In the sequel $h$ denotes the time step, and $t_n = n h$ are the
associated time instants, starting at $t_0=0$.
\section{Definition of $\RL_k$ schemes and consistency}
\label{sec:RL-def}
Let us consider a solution $y(t)$ of equation~\eqref{F2} on a time
interval $[0,T]$. It is recalled that the scheme~\eqref{eq:RL} is consistent
of order $k$ if, given a time step $h$, a time instant
$kh \le t_n \le T-h$, and the numerical approximation $y_{n+1}$ in~\eqref{eq:RL}
computed with $y_{n-j}=y(t_{n-j})$ for $j=0,\ldots, k-1$,
% \begin{itemize}
% \item \blue{giv en} a time step $h$ and a time instant $kh\le t_n \le T-h$,
% \item \blue{given} the numerical approximation $y_{n+1}$ in~\eqref{eq:RL}
%   computed with $y_{n-j}=y(t_{n-j})$ for $j=0,\ldots, k-1$,
% \end{itemize}
we have $|y_{n+1} - y(t_{n}+h)| \le C h^{k+1}$, for a constant $C$ only
depending on the data $a$, $b$, $y_0$ and $T$ of the problem~\eqref{F2}.
\begin{lemma}
  \label{prop:consistency}
  Assume that the functions $a(t,y)$ and $b(t,y)$ 
  are $\mathcal{C}^k$
  regular on $[0,T]\times \R^N$. 
  Moreover, assume that $a(t,y)$ is
  diagonal ($a(t,y) = \diag\left(a_i(t,y)\right)$) or constant. 
  Then the scheme in~\eqref{eq:RL} is consistent of order $k$
  for $k = 2$, $3$, $4$ if
  \begin{itemize}
  \item for $k=2$, we have
    \begin{equation*}
      \alpha_n = a_n + \frac{1}{2} a'_n h + O(h^2), \quad\text{and } \beta_n =
      b_n + \frac{1}{2} b'_n h + O(h^2);
    \end{equation*}
  \item for $k=3$, we have
    \begin{align*}
      \alpha_n = a_n + \frac{1}{2} a'_n h + \frac{1}{6} a''_n h^2 + O(h^3),
      \\ \text{and } \beta_n = b_n + \frac{1}{2} b'_n h +
      \frac{1}{12}(a'_nb_n-a_nb'_n) h^2 + O(h^3);
    \end{align*}
  \item for $k=4$, we have
    \begin{align*}
      \alpha_n = a_n + \frac{1}{2} a'_n h + \frac{1}{6} a''_n h^2 +
      \frac{1}{24}a'''_n h^3 + O(h^4), \\ \text{and } 
      \beta_n = b_n + \frac{1}{2} b'_n h + \frac{1}{12}(a'_nb_n-a_nb'_n) h^2 + \\
      \frac{1}{24} \left( b'''_n + a''_n b_n-a_n b''_n \right) h^3 + O(h^4); 
    \end{align*}
  \end{itemize}
  where $a'_n$, $a''_n$, $a'''_n$ and $b'_n$, $b''_n$, $b'''_n$ denote the
  successive derivatives at time $t_n$ of the functions $t\mapsto a(t,y(t))$ and
  $t\mapsto b(t,y(t))$.
\end{lemma}
\begin{remark}\label{rem:2}
  The assumption \textit{``$a(t,y)$ is diagonal or constant''} in
  Lemma~\ref{prop:consistency} has the following origin. To analyze the
  consistency of the scheme, we will compute a Taylor expansion
  in $h$ of the scheme in~\eqref{eq:RL}. This expansion is
  derived from Taylor expansions of $\alpha_n$ and $\beta_n$. For the
  sake of simplicity, assume the simple form
  $\alpha_n = \alpha_{n,0} + h\alpha_{n,1}$. We need to expand
  $\phi_1(\alpha_n h)$ as a power series in $h$, where the function
  $\phi_1$ is analytic on $\C$. However, in the matrix case, the
  equality,
  $\phi_1(M+N) = \phi_1(M) + \phi_1'(M)N + \ldots + \phi_1^{(i)}(M)N^i / i! +
  \ldots $ holds if $M$ and $N$ are commutative matrices.  Therefore one cannot
  expand $\phi_1(\alpha_n h)$ without the assumptions that $\alpha_{n,0}$ and
  $\alpha_{n,1}$ commute. This difficulty vanishes if $a(t,y)$ is
  constant or a varying diagonal matrix. 
  % We finally point out that, because
  % of the target application in Section \ref{sec:num-res}, we are specifically
  % interested \blue{in the case of a diagonal, non-constant, matrix $a(t,y)$}.
  % \blue{Therefore we did not perform numerical tests in
  % the case where $a(t,y)$ is constant.}
  % which case has been considered because it fits in the theoretical study.}
  % we did not perform numerical tests in that case. It involves a technical
  % difficulty mentioned in the introduction: the computation of
  % $\phi_1(\alpha_n h)v$.
  % That difficulty is general for exponential integrators and has been widely
  % studied,
  % we refer to \cite[Section 4]{Hochbruck-Ostermann-2010} on that point.
\end{remark}
\begin{proof}
  Consider equation~\eqref{F2} on the closed time interval $[0,T]$, and
  its solution, the function $y$. Since the functions $a$ and $b$ are
  $\mathcal{C}^k$ regular on $[0,T]\times \R^N$, the solution $y(t)$ is
  $\mathcal{C}^{k+1}$ regular on $[0,T]$. Its derivatives up to order
  $k+1$ are bounded by constants only depending on the data of 
  problem~\eqref{F2}, and on $T$. The Taylor expansion of $y$ at time instant $t_n$ is 
  \begin{equation*}
    % \label{taylor_ex_sol}
    y(t_n+h) = y(t_n)  + \sum_{j=1}^k \frac{s_j}{j!} h^j  + O(h^{k+1}),
  \end{equation*}
  with $s_j=y^{(j)}(t_n)$. Using that $y'=ay+b$ we get that
  \begin{align*}
    s_1 = &a_n y_n + b_n, \\
    s_2 = &(a'_n + a_n^2)y_n + a_nb_n + b'_n , \\
    s_3=&(a''_n + 3 a_n a'_n + a_n^3 )y_n 
          + 
          b''_n + a_n b'_n + 2 a'_n b_n + a_n^2 b_n , \\
    s_4 = &(a'''_n + 4 a''_n a_n + 3a_n^{' 2} + 6a'_n a_n^2 + a_n^4)y_n \\
          & \mbox{}\hskip 3em + b'''_n + b''_n a_n  +
            3 a''_n  b_n + 5 a'_n a_n b_n + 3 a'_n b'_n + a_n^3 b_n +a^2_n b'_n   .
  \end{align*}
  Series expansions %A series expansion
  in $h$ for $\alpha_n$ and for $\beta_n$
  are %is
  introduced as
  % \blue{Consider an expansion of the sequences of coefficients $\alpha_n$ and
  % $\beta_n$ of the form}
  \begin{align*}
    \label{general_form}
    \alpha_n &= \alpha_{n,0} + \alpha_{n,1} h + \dots + \alpha_{n,k-1}h^{k-1} + O(h^k), \\
    \beta_n  &= \beta_{n,0} + \beta_{n,1} h  + \dots + \beta_{n,k-1}h^{k-1} + O(h^k).
  \end{align*}
  If the matrix $a(t,y)$ is diagonal or constant (see Remark
  \ref{rem:2}), the Taylor expansion of the numerical solution $y_{n+1}$ in~\eqref{eq:RL} can be performed
  \begin{equation*}
    \label{taylor_num_sol}
    y_{n+1} = y(t_n) + \sum_{j=1}^k \dfrac{r_j}{j!} h^j + O(h^{k+1}).
  \end{equation*}
  A direct computation of the $r_j$ gives
  \begin{align*}
    % \label{eq:def-rj}
    r_1 = &\alpha_{n,0} y_n + \beta_{n,0}, 
    \\[5pt] \notag 
    r_2 = &(2\alpha_{n,1} + \alpha_{n,0}^2)y_n + 2\beta_{n,1} + \alpha_{n,0} \beta_{n,0} ,
    \\[5pt] \notag
    r_3=&(6\alpha_{n,2} + \alpha_{n,0}^3 +6 \alpha_{n,0} \alpha_{n,1} )y_n 
          + 
          3\alpha_{n,1} \beta_{n,0} + 6\beta_{n,2} + \alpha_{n,0}^2 \beta_{n,0} +3 \alpha_{n,0} \beta_{n,1} ,
    \\[5pt] \notag
    r_4 = &(24 \alpha_{n,0} \alpha_{n,2} + 24 \alpha_{n,3} + 12 \alpha_{n,1} \alpha_{n,0}^2 + 12 \alpha_{n,1} ^2 + \alpha_{n,0}^4)y_n
    \\ &\notag
         +12 \alpha_{n,2} \beta_{n,0}   + 24 \beta_{n,3} + 12 \alpha_{n,0} \beta_{n,2} +12 \alpha_{n,1} \beta_{n,1} +4 \alpha_{n,0}^2 \beta_{n,1} 
         + 8 \alpha_{n,0} \alpha_{n,1} \beta_{n,0} + \alpha_{n,0}^3 \beta_{n,0},
  \end{align*}
  where $y_n$ denotes $y(t_n)$.
  The condition to be consistent of order $k$ is:
  $r_i =s_i$ for $1 \le i \le k$.
  The consistency conditions in
  Lemma \ref{prop:consistency}
  are obtained by
  solving recursively these relations.
\end{proof}

We then can state our main result, which includes the definition of the
$\RL_k$ schemes.
\begin{theorem}
  \label{def:rlk}
  Assume (as in Lemma \ref{prop:consistency}) that the functions $a(t,y)$
  and $b(t,y)$ are $\mathcal{C}^k$ regular 
  on $[0,T]\times \R^N$,
  and that $a(t,y)$ is diagonal or constant. Then, the three
  schemes defined for $k = 2$, $3$, $4$ by equation~\eqref{eq:RL} and the
  coefficients,
  \begin{itemize}
  \item for $k=2$,
    \begin{equation*}
      \alpha_n = \frac{3}{2}a_n-\frac{1}{2}a_{n-1}, 
      \beta_n = \frac{3}{2}b_n-\frac{1}{2}b_{n-1},
    \end{equation*}
  \item for $k=3$,
    \begin{align*}
      \alpha_n &= \frac{1}{12}(23a_n - 16 a_{n-1} + 5 a_{n-2}), 
      \\
      \beta_n &= \frac{1}{12}(23b_n - 16 b_{n-1} + 5 b_{n-2}) +
                \frac{h}{12}(a_nb_{n-1} - a_{n-1}b_n),
    \end{align*}
  \item for $k=4$,
    \begin{align*}
      \alpha_n &= \frac{1}{24}(55a_n - 59 a_{n-1} + 37 a_{n-2}-9a_{n-3}), \\
      \beta_n &= \frac{1}{24}(55b_n - 59 b_{n-1} + 37 b_{n-2}-9b_{n-3}) \\ 
               &+ \frac{h}{12}(a_n(3b_{n-1}-b_{n-2}) - (3a_{n-1}-a_{n-2})b_n),
    \end{align*}
  \end{itemize}
  where $a_{j}=a(t_{j}, y_{j})$ and $b_{j}=b(t_{j}, y_{j})$, are consistent of
  order $k$.
  \\
  The three methods stated above are called Rush-Larsen methods of order
  $k$, and denoted by $\RL_k$. They are explicit and $k$-step methods.
\end{theorem}
\begin{remark}
  If the matrix $a$ is a constant, $a(t,y)=A$, then we have
  $\alpha_n = A$ for all three methods. In this case, the
  expressions of the coefficients $\beta_n$ in Theorem~\ref{def:rlk} for
  $k=3$, $4$ simplify as follows:
  \begin{align*}
    \RL_3\text{~~case}:\qquad  
    \beta_n &= \dfrac{1}{12}(23b_n - 16 b_{n-1} + 5 b_{n-2}) - \dfrac{h}{12}A(b_n - b_{n-1}).
    \\
    \RL_4 \text{~~case}:\qquad  
    \beta_n &= \dfrac{1}{24}(55b_n - 59 b_{n-1} + 37 b_{n-2}-9b_{n-3}) 
              - \dfrac{h}{12}A(2b_n - 3b_{n-1} + b_{n-2}).
  \end{align*}
\end{remark}
\begin{proof}
  It is a direct consequence of backwards differentiation formulas, that
  we first recall. The derivatives of a real function $f$ at the
  time instant $t_n$ can be approximated as follows (with common notations):
  \begin{itemize}
  \item first derivative,
    \begin{align*}
      f'_n & = \frac{f_n-f_{n-1}}{h} + O(h) \\
           & = \frac{1}{2h}\left( 3f_n - 4 f_{n-1} + f_{n-2} \right) + O(h^2) \\
           & = \frac{1}{6h}\left( 11f_n - 18 f_{n-1} + 9f_{n-2} - 2 f_{n-3} \right) + O(h^3);
    \end{align*}
  \item second derivative,
    \begin{align*}
      f''_n & = \frac{1}{h^2}\left( f_n -2f_{n-1} +f_{n-2} \right) + O(h) \\
            & = \frac{1}{h^2}\left( 2f_n - 5 f_{n-1} + 4f_{n-2} -  f_{n-3} \right) + O(h^2);
    \end{align*}
  \item third derivative,
    \begin{equation*}
      f'''_n  = \frac{1}{h^3}\left( f_n - 3 f_{n-1} + 3f_{n-2} -  f_{n-3} \right) + O(h).   
    \end{equation*}
  \end{itemize}
  With these formulas, the consistency condition at order 3 on
  the coefficient $\alpha_n$ becomes
  \begin{align*}
    \alpha_n  &= a_n +  \frac{1}{2} a'_n h + \frac{1}{6} a''_n h^2 + O(h^3) \\
              &= a_n +  \frac{1}{4} \left( 3a_n - 4 a_{n-1} + a_{n-2} \right)
                + \frac{1}{6}  \left( a_n -2a_{n-1} +a_{n-2} \right) + O(h^3) \\
              &= \frac{1}{12}(23a_n - 16 a_{n-1} + 5 a_{n-2})  + O(h^3).     
  \end{align*}
  We retrieve the definition of $\alpha_n$ for the $\RL_3$ scheme. The same
  proof holds for $\beta_n$, and extends to order 4.
\end{proof}

\section{Stability under perturbation and convergence}
\label{sec:conv}
We refer to \cite[Ch. III-8]{Hairer-ODE-I} for the definitions of convergence
and of stability under perturbation.  For the analysis of time-stepping methods,
it is commonly assumed that $f$ in equation~\eqref{F1} is uniformly
Lipschitz with respect to its second variable $y$. This
hypothesis will be replaced by assumptions based on the formulation~\eqref{F2}.
Precisely it will be assumed that
\begin{equation}
  \label{eq:a-b-f-Lipschitz}
  a(t,y) ~~\text{is bounded},\quad 
  a(t,y),~b(t,y)~~\text{are uniformly Lipschitz in}~y.
\end{equation}
The Lipschitz constants of $a$ and $b$ are denoted by $L_a$ and
$L_b$, respectively. The upper bound on $|a(t,y)|$ is denoted by $M_a$.
\begin{proposition}
  \label{prop:stab}
  If the assumption~\eqref{eq:a-b-f-Lipschitz} holds, then, the $\RL_k$
  schemes are stable under perturbation, for $k=2$, $3$, $4$.
  In addition, also for $k=2$, $3$, $4$, if the consistency assumptions of Theorem~\ref{def:rlk}
  are satisfied ($a(t,y)$ and $b(t,y)$ are $\mathcal{C}^k$ regular and $a(t,y)$
  is diagonal or constant), then the $\RL_k$ scheme is convergent of
  order $k$.
\end{proposition}

Stability under perturbation together with consistency implies (nonstiff)
convergence, see \textit{e.g.} \cite{Hairer-ODE-I}, or
\cite{coudiere-lontsi-pierre-2015}, where the current setting has been detailed.
Therefore the proof of the convergence statement in
Proposition~\ref{prop:stab} 
is immediate, and will not be recalled here.

The following definitions are necessary to prove Proposition \ref{prop:stab}.
Equation~\eqref{F1} is considered on $E = \R^N$ with the max norm $|\cdot|$. A
final time $T>0$ is considered.  The space of $N \times N$ matrices is equipped
with the operator norm $\Vert\cdot\Vert$ associated with $|\cdot|$.
The space $E^k$ is equipped with the max norm
$|Y|_\infty = \max_{1\le i\le k} |y_i|$ with $Y=(y_1,\ldots,y_k)$. The $\RL_k$
scheme is defined by the mapping
\begin{displaymath}
  s\th:~ Y=(y_1,\ldots,y_k) \in E^k \longrightarrow
  s\th(Y) \in E,
\end{displaymath}
with
\begin{displaymath}
  s\th(Y)
  = y_k + h \phi_1(\alpha_{t,h}(Y)h)\left (
    \alpha\th(Y) y_k + \beta\th(Y)\right ), 
\end{displaymath}
in such a way that the scheme in~\eqref{eq:RL} reads
$y_{n+1} = s_{t_n,h}(y_{n-k+1}, \ldots, y_n)$. The functions $\alpha\th$ and
$\beta\th$ map the vector $Y$ of the $k$ previous values to the values
$\alpha_n$ and $\beta_n$ given in Theorem~\ref{def:rlk}. For instance,
the function $\alpha\th(Y)$ for $k=3$ (the $\RL_3$ scheme) reads
\begin{equation*}
  \alpha\th(Y) = \frac{1}{12}(23 a(t,y_3) - 16 a(t-h,y_2) + 5 a(t-2h, y_1)).
\end{equation*}  
% The scheme generator is the mapping $S\th$ given by, 
% \begin{displaymath}
%   S\th~:  Y=(y_1,\ldots,y_k) \in E^k \longrightarrow 
%   \left( y_2 , \ldots , y_k, s\th(Y)\right) \in    E^k .
% \end{displaymath}

A first technique to prove the stability under perturbation
consists in showing that the function $s\th$ is globally
Lipschitz in $Y$. To this aim, the derivative $\partial_Y s\th$ has to
be analyzed. As developed in Remark~\ref{rem:2}, it implies restrictions on
the function $a(t,y)$: it has to be either diagonal or
constant. A second technique consists in proving the following two
stability conditions:
\begin{gather}
  \label{eq:cond-stab-1}
  | s\th(Y) - s\th(Z) | \le 
  | Y - Z|_\infty
  \left( 1+C h (| Y |_\infty + 1 ) \right),
  \\
  \label{eq:cond-stab-2}
  | s\th(Y) | \le | Y |_\infty (1+C h)  + C h,
\end{gather}
for all $Y$ and $Z$ in $E^k$, and where the constant $C$ depends only on
the data $a$, $b$, $y_0$ in equation~\eqref{F2}, and on the final time $T$.
These are sufficient conditions for the stability under perturbation, as proved
in \cite[Section 2]{coudiere-lontsi-pierre-2015}.
We will use the conditions~\eqref{eq:cond-stab-1}
and~\eqref{eq:cond-stab-2} here, because they are more general, and give rise
to less computations. The core of the proof is the following property of the
$\RL_k$ scheme.  For $Y=(y_1,\ldots,y_k) \in E^k$, we have
\begin{equation}
  \label{eq:rlk-diff-form}
  s\th(Y) = z(t+h) \quad  \text{for }
  z' = \alpha\th(Y) z + \beta\th(Y), \quad  z(t)=y_k.
\end{equation}
It will be used together with the following Gronwall inequality (see \cite[Lemma
196, p.150]{gronwall-monograph}). Suppose that $z(t)$ is a $\mathcal{C}^1$
function, and that there exists $M_1>0$ and $M_2>0$ such that
$|z'(t)| \le M_1 |z(t)| + M_2$ for all $t\in[t_0,t_0+h]$. Then
\begin{equation}
  \label{eq:gronwall}
  \forall t \in [t_0,t_0+h],\quad |z(t)| \le \e^{M_1 (t-t_0)} \left( |z(t_0)| +
    M_2 (t-t_0) \right).
\end{equation}
\begin{proof}[Proposition \ref{prop:stab}]
  In this proof, we always assume that $0\le h,t\le T$, and denote by $C_i$ a
  constant that depends only on the data $a$, $b$ and $T$ of problem~\eqref{F2}.
  With the assumptions in~\eqref{eq:a-b-f-Lipschitz}, and the
  definitions of $\alpha_n$ ($k=2$, $3$, $4$) in Theorem~\ref{def:rlk}, the
  function $\alpha\th$ is uniformly Lipschitz with a Lipschitz constant equal to
  $L_\alpha$. Moreover we have the uniform bound
  $\Vert \alpha\th\Vert \le M_\alpha$. Since the function $b(t,y)$ is
  uniformly Lipschitz with respect to $y$, and since $0\le t \le T$, we
  have
  \begin{equation}
    \label{eq:Kb}
    | b(t,y)| \le |b(t,0)| + |b(t,y) - b(t,0)| \le K_b( 1 + | y| ),
  \end{equation}
  with $K_b = \max(L_b, \sup_{0\le t\le T} |b(t,0)|)$. For the $\RL_3$
  scheme, we have
  \begin{equation*}
    |\beta\th(Y)|_\infty \le \frac{11}{3} K_b( 1 + | Y|_\infty )
    + \frac{h}{12} M_a 2 K_b( 1 + | Y|_\infty ) \le C_1 ( 1 + | Y|_\infty ).
  \end{equation*}
  The same inequality holds for the $\RL_2$, 
  and $\RL_4$ schemes. 
  Afterwards, we can apply these bounds to
  the differential equation in~\eqref{eq:rlk-diff-form}
  $$| z ' | = | \alpha\th(Y)z + \beta\th(Y)| \le M_\alpha |z| + C_1( 1 + |
  Y|_\infty ).$$ The initial state is $|z(t)| = |y_k| \le
  |Y|_\infty$. Finally, the Gronwall inequality~\eqref{eq:gronwall}
  yields, for $t\le \tau \le t+h$,
  \begin{align}
    \label{ineq1-proof-stab}
    |z(\tau)| &\le \e^{M_\alpha h} \left ( |Y|_\infty + h C_1( 1 + | Y|_\infty
                )\right ) 
                \notag \\ \notag
              &\le \e^{M_\alpha h} \left ( |Y|_\infty (1 + C_1h) + C_1h \right
                ) \\ &\le |Y|_\infty (1 + C_2 h) + C_2 h,
  \end{align}
  by bounding the exponential with an affine function for $0\le h\le T$. This
  gives the stability condition~\eqref{eq:cond-stab-2} for $\tau = t+h$.

  For the $\RL_2$ scheme, the function $\beta\th$ is uniformly
  Lipschitz. For the $\RL_3$ scheme, for $Y = (y_1,y_2,y_3)$ and
  $Z = (z_1,z_2,z_3)$ in $E^3$, we have
  \begin{align*}
    \vert\beta\th(Y)- \beta\th(Z)\vert_\infty  \le
    \dfrac{11}{3}L_b\vert Y-Z\vert_\infty  + \dfrac{h}{12}
    \big( 
    &\left \vert 
      a(t,y_3)b(t-h,y_2)
      -
      a(t,z_3)b(t-h,z_2)
      \right \vert 
      \big.
    \\+
    \big .
    &\left \vert 
      a(t-h,y_2)b(t,y_3)
      -
      a(t-h,z_2)b(t,z_3)
      \right \vert 
      \big)
  \end{align*}
  Let us bound the Lipschitz constant of a function of the type
  $F(Y)=a(\xi,y_2)b(\tau,y_3)$, for $0\le \tau, \xi\le T$:
  \begin{align*}
    \left \vert 
    F(Y) - F(Z)
    \right \vert &= 
                   \left \vert 
                   a(\xi,y_3)\left (
                   b(\tau,y_2) - b(\tau,z_2)
                   \right )
                   +
                   \left (
                   a(\xi,y_3)-a(\xi,z_3)
                   \right )
                   b(\tau,z_2)
                   \right \vert 
    \\ &\le
         M_a L_b \vert Y-Z\vert _\infty 
         +
         L_a \vert Y-Z\vert _\infty \vert b(\tau,z_2)\vert .
  \end{align*}
  With the inequality~\eqref{eq:Kb}, this yields, for
  $0\le \tau,\xi\le T$, and $Y$, $Z$ in $E^k$,
  $\left| F(Y) - F(Z) \right | \le C_3 | Y-Z|_\infty (1 + | Z| _\infty)$. As a
  result, we have
  \begin{equation*}
    |\beta\th(Y)- \beta\th(Z)|_\infty \le C_4 |Y-Z|_\infty \left( 1 + |
      Z|_\infty \right).
  \end{equation*}
  The same inequality holds for the $\RL_4$ scheme.

  Finally we consider $Y_1$ and $Y_2$ in $E^k$, and the notation
  $\alpha_i = \alpha\th(Y_i)$, and $\beta_i = \beta\th(Y_i)$. The
  property~\eqref{eq:rlk-diff-form} shows that
  $s\th(Y_1)-s\th(Y_2) = (z_1-z_2)(t+h)$, where $z_i$ is the solution to
  $z'_i = \alpha_i z_i + \beta_i$, with $z_i(t) = Y_{i,k}$. On the first hand,
  with the inequality~\eqref{ineq1-proof-stab}, we have
  $|z_2(\tau)| \le C_5(1+|Y_2|_\infty )$ for $t\le\tau\le t+h$. On the second
  hand, on $[t,t+h]$, we have
  \begin{align*}
    |(z_1-z_2)'| &\le
                   |\alpha_1| |z_1-z_2| + |\alpha_1-\alpha_2| |z_2|
                   + |\beta_1-\beta_2|
    \\ &\le
         M_\alpha |z_1-z_2| + L_\alpha |Y_1-Y_2|_\infty  C_5(1+|Y_2|_\infty ) + C_4 |Y_1-Y_2|_\infty  (1+|Y_2|_\infty )
    \\ &\le M_\alpha |z_1-z_2| + C_6|Y_1-Y_2|_\infty  (1+|Y_2|_\infty ).
  \end{align*}
  The initial condition yields
  $|(z_1-z_2)(t)| = |Y_{1,k}-Y_{2,k}|\le |Y_1-Y_2|_\infty$. As a
  consequence, the Gronwall inequality~\eqref{eq:gronwall} applied to
  these bounds shows that
  \begin{align*}
    |(z_1-z_2)(t+h)| \le &\e^{M_\alpha h}
                           \left (
                           |Y_1-Y_2|_\infty + h C_6|Y_1-Y_2|_\infty  (1+|Y_2|_\infty )
                           \right ) 
    \\ \le&
            \e^{M_\alpha h} |Y_1-Y_2|_\infty 
            \left (
            1+C_6 h (1+|Y_2|_\infty
            ) \right ).
  \end{align*}
  This last inequality implies the stability condition~\eqref{eq:cond-stab-1},
  again by bounding the exponential with an affine function for $0\le h\le T$.
\end{proof}
\section{Dahlquist stability}
\label{sec:dahlquist-stab}
For the general ideas and definitions concerning the Dahlquist stability we
refer to \cite{Hairer-ODE-II}.  The background for the Dahlquist stability of
exponential integrators with a general varying stabilizer $a(t,y)$ has been
developed in \cite{coudiere-lontsi-pierre-2015}, following the ideas of Perego
and Veneziani \cite{perego-2009}. The equation~\eqref{F1} is considered
with the Dahlquist test function $f(t,y)=\lambda y$, which is split into
$f(t,y)= a(t,y)y + b(t,y)$, in order to match the framework of
equation~\eqref{F2}, with
\begin{equation*}
  a(t,y) = \theta \lambda, \quad b(t,y) = \lambda(1-\theta)y.
\end{equation*}
For $\theta = 1$, the methods are exact and thus $A$-stable.
For $\theta\simeq 1$, the exact linear part of $f(t,y)$ in
equation~\eqref{F1} is well approximated by $a(t,y)$. The stability domain
depends on $\theta$, it is denoted by $D_\theta$. Given a value
of $\theta$, the region $D_\theta$ is defined by the modulus of a
stability function, with the same definition as for multistep methods, see
\textit{e.g.} \cite{Hairer-ODE-II}. This stability function has been numerically computed, pointwise on a grid inside the complex plane $\C$,
for each of the three $\RL_k$ schemes, $k=2$, $3$, $4$.

\subsubsection*{Order 2 Rush-Larsen}
\begin{figure}
  \centering
  % GNUPLOT: LaTeX picture with Postscript
  \begingroup
  \sizeg
  \makeatletter
  \providecommand\color[2][]{%
    \GenericError{(gnuplot) \space\space\space\@spaces}{%
      Package color not loaded in conjunction with
      terminal option `colourtext'%
    }{See the gnuplot documentation for explanation.%
    }{Either use 'blacktext' in gnuplot or load the package
      color.sty in LaTeX.}%
    \renewcommand\color[2][]{}%
  }%
  \providecommand\includegraphics[2][]{%
    \GenericError{(gnuplot) \space\space\space\@spaces}{%
      Package graphicx or graphics not loaded%
    }{See the gnuplot documentation for explanation.%
    }{The gnuplot epslatex terminal needs graphicx.sty or graphics.sty.}%
    \renewcommand\includegraphics[2][]{}%
  }%
  \providecommand\rotatebox[2]{#2}%
  \@ifundefined{ifGPcolor}{%
    \newif\ifGPcolor
    \GPcolortrue
  }{}%
  \@ifundefined{ifGPblacktext}{%
    \newif\ifGPblacktext
    \GPblacktexttrue
  }{}%
  % define a \g@addto@macro without @ in the name:
  \let\gplgaddtomacro\g@addto@macro
  % define empty templates for all commands taking text:
  \gdef\gplbacktext{}%
  \gdef\gplfronttext{}%
  \makeatother
  \ifGPblacktext
  % no textcolor at all
  \def\colorrgb#1{}%
  \def\colorgray#1{}%
  \else
  % gray or color?
  \ifGPcolor
  \def\colorrgb#1{\color[rgb]{#1}}%
  \def\colorgray#1{\color[gray]{#1}}%
  \expandafter\def\csname LTw\endcsname{\color{white}}%
  \expandafter\def\csname LTb\endcsname{\color{black}}%
  \expandafter\def\csname LTa\endcsname{\color{black}}%
  \expandafter\def\csname LT0\endcsname{\color[rgb]{1,0,0}}%
  \expandafter\def\csname LT1\endcsname{\color[rgb]{0,1,0}}%
  \expandafter\def\csname LT2\endcsname{\color[rgb]{0,0,1}}%
  \expandafter\def\csname LT3\endcsname{\color[rgb]{1,0,1}}%
  \expandafter\def\csname LT4\endcsname{\color[rgb]{0,1,1}}%
  \expandafter\def\csname LT5\endcsname{\color[rgb]{1,1,0}}%
  \expandafter\def\csname LT6\endcsname{\color[rgb]{0,0,0}}%
  \expandafter\def\csname LT7\endcsname{\color[rgb]{1,0.3,0}}%
  \expandafter\def\csname LT8\endcsname{\color[rgb]{0.5,0.5,0.5}}%
  \else
  % gray
  \def\colorrgb#1{\color{black}}%
  \def\colorgray#1{\color[gray]{#1}}%
  \expandafter\def\csname LTw\endcsname{\color{white}}%
  \expandafter\def\csname LTb\endcsname{\color{black}}%
  \expandafter\def\csname LTa\endcsname{\color{black}}%
  \expandafter\def\csname LT0\endcsname{\color{black}}%
  \expandafter\def\csname LT1\endcsname{\color{black}}%
  \expandafter\def\csname LT2\endcsname{\color{black}}%
  \expandafter\def\csname LT3\endcsname{\color{black}}%
  \expandafter\def\csname LT4\endcsname{\color{black}}%
  \expandafter\def\csname LT5\endcsname{\color{black}}%
  \expandafter\def\csname LT6\endcsname{\color{black}}%
  \expandafter\def\csname LT7\endcsname{\color{black}}%
  \expandafter\def\csname LT8\endcsname{\color{black}}%
  \fi
  \fi
  \setlength{\unitlength}{0.0500bp}%
  \begin{picture}(9070.00,3968.00)%
    \gplgaddtomacro\gplbacktext{%
      \csname LTb\endcsname%
      \put(396,440){\makebox(0,0)[r]{\strut{} 0}}%
      \put(396,826){\makebox(0,0)[r]{\strut{} 1}}%
      \put(396,1212){\makebox(0,0)[r]{\strut{} 2}}%
      \put(396,1598){\makebox(0,0)[r]{\strut{} 3}}%
      \put(396,1984){\makebox(0,0)[r]{\strut{} 4}}%
      \put(396,2369){\makebox(0,0)[r]{\strut{} 5}}%
      \put(396,2755){\makebox(0,0)[r]{\strut{} 6}}%
      \put(396,3141){\makebox(0,0)[r]{\strut{} 7}}%
      \put(396,3527){\makebox(0,0)[r]{\strut{} 8}}%
      \put(528,220){\makebox(0,0){\strut{}-6}}%
      \put(1428,220){\makebox(0,0){\strut{}-5}}%
      \put(2327,220){\makebox(0,0){\strut{}-4}}%
      \put(3227,220){\makebox(0,0){\strut{}-3}}%
      \put(4126,220){\makebox(0,0){\strut{}-2}}%
      \put(5026,220){\makebox(0,0){\strut{}-1}}%
      \put(5925,220){\makebox(0,0){\strut{} 0}}%
      \put(6825,220){\makebox(0,0){\strut{} 1}}%
    }%
    \gplgaddtomacro\gplfronttext{%
    }%
    \gplgaddtomacro\gplbacktext{%
      \csname LTb\endcsname%
      \put(396,440){\makebox(0,0)[r]{\strut{} 0}}%
      \put(396,826){\makebox(0,0)[r]{\strut{} 1}}%
      \put(396,1212){\makebox(0,0)[r]{\strut{} 2}}%
      \put(396,1598){\makebox(0,0)[r]{\strut{} 3}}%
      \put(396,1984){\makebox(0,0)[r]{\strut{} 4}}%
      \put(396,2369){\makebox(0,0)[r]{\strut{} 5}}%
      \put(396,2755){\makebox(0,0)[r]{\strut{} 6}}%
      \put(396,3141){\makebox(0,0)[r]{\strut{} 7}}%
      \put(396,3527){\makebox(0,0)[r]{\strut{} 8}}%
      \put(528,220){\makebox(0,0){\strut{}-6}}%
      \put(1428,220){\makebox(0,0){\strut{}-5}}%
      \put(2327,220){\makebox(0,0){\strut{}-4}}%
      \put(3227,220){\makebox(0,0){\strut{}-3}}%
      \put(4126,220){\makebox(0,0){\strut{}-2}}%
      \put(5026,220){\makebox(0,0){\strut{}-1}}%
      \put(5925,220){\makebox(0,0){\strut{} 0}}%
      \put(6825,220){\makebox(0,0){\strut{} 1}}%
    }%
    \gplgaddtomacro\gplfronttext{%
      \csname LTb\endcsname%
      \put(7158,3390){\makebox(0,0)[l]{\strut{}$\theta$ = 0}}%
      \csname LTb\endcsname%
      \put(7158,3115){\makebox(0,0)[l]{\strut{}$\theta$ = 0.5}}%
      \csname LTb\endcsname%
      \put(7158,2840){\makebox(0,0)[l]{\strut{}$\theta$ = 2/3}}%
      \csname LTb\endcsname%
      \put(7158,2565){\makebox(0,0)[l]{\strut{}$\theta$ = 0.7}}%
      \csname LTb\endcsname%
      \put(7158,2290){\makebox(0,0)[l]{\strut{}$\theta$ = 5/6}}%
      \csname LTb\endcsname%
      \put(7158,2015){\makebox(0,0)[l]{\strut{}$\theta$ = 2}}%
    }%
    \gplgaddtomacro\gplbacktext{%
      \csname LTb\endcsname%
      \put(396,440){\makebox(0,0)[r]{\strut{} 0}}%
      \put(396,826){\makebox(0,0)[r]{\strut{} 1}}%
      \put(396,1212){\makebox(0,0)[r]{\strut{} 2}}%
      \put(396,1598){\makebox(0,0)[r]{\strut{} 3}}%
      \put(396,1984){\makebox(0,0)[r]{\strut{} 4}}%
      \put(396,2369){\makebox(0,0)[r]{\strut{} 5}}%
      \put(396,2755){\makebox(0,0)[r]{\strut{} 6}}%
      \put(396,3141){\makebox(0,0)[r]{\strut{} 7}}%
      \put(396,3527){\makebox(0,0)[r]{\strut{} 8}}%
      \put(528,220){\makebox(0,0){\strut{}-6}}%
      \put(1428,220){\makebox(0,0){\strut{}-5}}%
      \put(2327,220){\makebox(0,0){\strut{}-4}}%
      \put(3227,220){\makebox(0,0){\strut{}-3}}%
      \put(4126,220){\makebox(0,0){\strut{}-2}}%
      \put(5026,220){\makebox(0,0){\strut{}-1}}%
      \put(5925,220){\makebox(0,0){\strut{} 0}}%
      \put(6825,220){\makebox(0,0){\strut{} 1}}%
    }%
    \gplgaddtomacro\gplfronttext{%
    }%
    \gplbacktext
    \put(0,0){\includegraphics{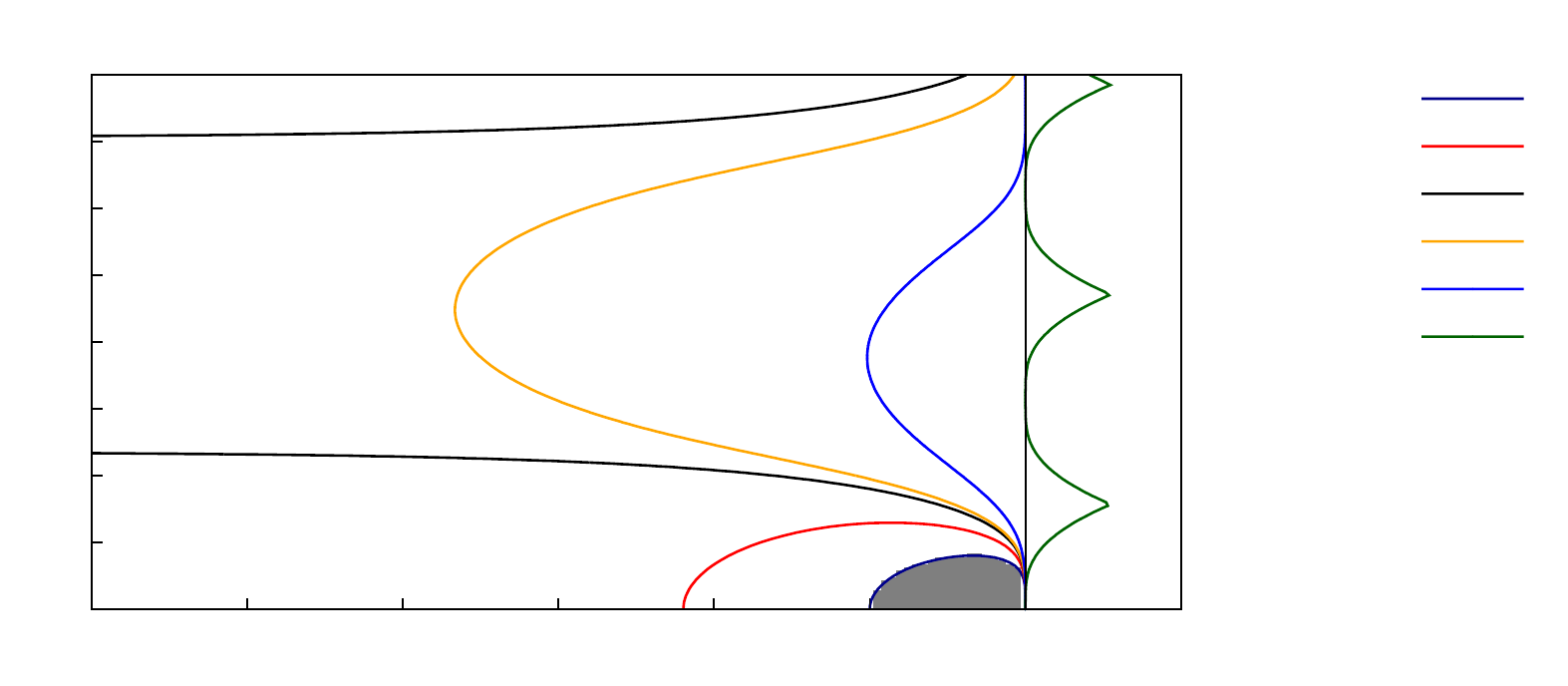}}%
    \gplfronttext
  \end{picture}%
  \endgroup
  \caption{Stability domain $D_\theta$ for the $\RL_2$ scheme for various values
    of $\theta$. The stability domain for the particular case $\theta=0$ (no
    stabilization) is in gray, corresponding to the Adams-Bashforth scheme of
    order 2. }
  \label{fig:RL2-D-stab}
\end{figure}

The stability domain for the $\RL_2$ scheme has been analyzed
in~\cite{perego-2009}. The situation for this scheme is interesting, and we
reproduced the results on Figure~\ref{fig:RL2-D-stab}. We note the
observations below.
\begin{itemize}
\item If $0\le \theta<2/3$, the stability domain $D_\theta$ is bounded. Its size
  increases with $\theta$, starting from the stability domain without
  stabilization for $\theta=0$, that corresponds to the Adams-Bashforth scheme
  of order 2.
\item If $\theta=2/3$, the method is $A(0)$ stable: $\R^- \subset D_\theta$.
  The domain boundary is asymptotically parallel to the real axis, so that the
  method is not $A(\alpha)$ stable.
\item If $\theta>2/3$, the stability domain is located around the $y$-axis: the
  method is $A(\alpha)$ stable. The angle $\alpha$ increases with $\theta$, it
  goes to $\pi/2$ as $\theta\to 1^-$.
\end{itemize}
\subsubsection*{Rush-Larsen methods of orders 3 and 4}
\begin{figure}
  \centering
  % GNUPLOT: LaTeX picture with Postscript
  \begingroup
  \sizeg
  \makeatletter
  \providecommand\color[2][]{%
    \GenericError{(gnuplot) \space\space\space\@spaces}{%
      Package color not loaded in conjunction with
      terminal option `colourtext'%
    }{See the gnuplot documentation for explanation.%
    }{Either use 'blacktext' in gnuplot or load the package
      color.sty in LaTeX.}%
    \renewcommand\color[2][]{}%
  }%
  \providecommand\includegraphics[2][]{%
    \GenericError{(gnuplot) \space\space\space\@spaces}{%
      Package graphicx or graphics not loaded%
    }{See the gnuplot documentation for explanation.%
    }{The gnuplot epslatex terminal needs graphicx.sty or graphics.sty.}%
    \renewcommand\includegraphics[2][]{}%
  }%
  \providecommand\rotatebox[2]{#2}%
  \@ifundefined{ifGPcolor}{%
    \newif\ifGPcolor
    \GPcolortrue
  }{}%
  \@ifundefined{ifGPblacktext}{%
    \newif\ifGPblacktext
    \GPblacktexttrue
  }{}%
  % define a \g@addto@macro without @ in the name:
  \let\gplgaddtomacro\g@addto@macro
  % define empty templates for all commands taking text:
  \gdef\gplbacktext{}%
  \gdef\gplfronttext{}%
  \makeatother
  \ifGPblacktext
  % no textcolor at all
  \def\colorrgb#1{}%
  \def\colorgray#1{}%
  \else
  % gray or color?
  \ifGPcolor
  \def\colorrgb#1{\color[rgb]{#1}}%
  \def\colorgray#1{\color[gray]{#1}}%
  \expandafter\def\csname LTw\endcsname{\color{white}}%
  \expandafter\def\csname LTb\endcsname{\color{black}}%
  \expandafter\def\csname LTa\endcsname{\color{black}}%
  \expandafter\def\csname LT0\endcsname{\color[rgb]{1,0,0}}%
  \expandafter\def\csname LT1\endcsname{\color[rgb]{0,1,0}}%
  \expandafter\def\csname LT2\endcsname{\color[rgb]{0,0,1}}%
  \expandafter\def\csname LT3\endcsname{\color[rgb]{1,0,1}}%
  \expandafter\def\csname LT4\endcsname{\color[rgb]{0,1,1}}%
  \expandafter\def\csname LT5\endcsname{\color[rgb]{1,1,0}}%
  \expandafter\def\csname LT6\endcsname{\color[rgb]{0,0,0}}%
  \expandafter\def\csname LT7\endcsname{\color[rgb]{1,0.3,0}}%
  \expandafter\def\csname LT8\endcsname{\color[rgb]{0.5,0.5,0.5}}%
  \else
  % gray
  \def\colorrgb#1{\color{black}}%
  \def\colorgray#1{\color[gray]{#1}}%
  \expandafter\def\csname LTw\endcsname{\color{white}}%
  \expandafter\def\csname LTb\endcsname{\color{black}}%
  \expandafter\def\csname LTa\endcsname{\color{black}}%
  \expandafter\def\csname LT0\endcsname{\color{black}}%
  \expandafter\def\csname LT1\endcsname{\color{black}}%
  \expandafter\def\csname LT2\endcsname{\color{black}}%
  \expandafter\def\csname LT3\endcsname{\color{black}}%
  \expandafter\def\csname LT4\endcsname{\color{black}}%
  \expandafter\def\csname LT5\endcsname{\color{black}}%
  \expandafter\def\csname LT6\endcsname{\color{black}}%
  \expandafter\def\csname LT7\endcsname{\color{black}}%
  \expandafter\def\csname LT8\endcsname{\color{black}}%
  \fi
  \fi
  \setlength{\unitlength}{0.0500bp}%
  \begin{picture}(9070.00,3968.00)%
    \gplgaddtomacro\gplbacktext{%
      \csname LTb\endcsname%
      \put(396,440){\makebox(0,0)[r]{\strut{} 0}}%
      \put(396,852){\makebox(0,0)[r]{\strut{} 20}}%
      \put(396,1263){\makebox(0,0)[r]{\strut{} 40}}%
      \put(396,1675){\makebox(0,0)[r]{\strut{} 60}}%
      \put(396,2086){\makebox(0,0)[r]{\strut{} 80}}%
      \put(396,2498){\makebox(0,0)[r]{\strut{} 100}}%
      \put(396,2910){\makebox(0,0)[r]{\strut{} 120}}%
      \put(396,3321){\makebox(0,0)[r]{\strut{} 140}}%
      \put(528,220){\makebox(0,0){\strut{}-250}}%
      \put(1787,220){\makebox(0,0){\strut{}-200}}%
      \put(3047,220){\makebox(0,0){\strut{}-150}}%
      \put(4306,220){\makebox(0,0){\strut{}-100}}%
      \put(5566,220){\makebox(0,0){\strut{}-50}}%
      \put(6825,220){\makebox(0,0){\strut{} 0}}%
      \put(6817,440){\rotatebox{240}{\makebox(0,0)[l]{\strut{}\textcolor{NavyBlue} {$\leftarrow$ -0.54 }}}}%
      \put(6470,440){\rotatebox{100}{\makebox(0,0)[l]{\strut{}\textcolor{red} {$\leftarrow$ -14.1 }}}}%
      \put(5928,440){\rotatebox{100}{\makebox(0,0)[l]{\strut{}\textcolor{black}{$\leftarrow$ -35.6}}}}%
      \put(4387,440){\rotatebox{100}{\makebox(0,0)[l]{\strut{}\textcolor{Orange}{$\leftarrow$ -96.8}}}}%
      \put(1042,440){\rotatebox{100}{\makebox(0,0)[l]{\strut{}\textcolor{blue}{$\leftarrow$ -229.6}}}}%
      \put(4049,440){\rotatebox{100}{\makebox(0,0)[l]{\strut{}\textcolor{OliveGreen}{$\leftarrow$ -110.2}}}}%
    }%
    \gplgaddtomacro\gplfronttext{%
      \csname LTb\endcsname%
      \put(7026,3390){\makebox(0,0)[l]{\strut{}$\theta$ = 0.85}}%
      \csname LTb\endcsname%
      \put(7026,3115){\makebox(0,0)[l]{\strut{}$\theta$ = 0.9}}%
      \csname LTb\endcsname%
      \put(7026,2840){\makebox(0,0)[l]{\strut{}$\theta$ = 0.95}}%
      \csname LTb\endcsname%
      \put(7026,2565){\makebox(0,0)[l]{\strut{}$\theta$ = 1.05}}%
      \csname LTb\endcsname%
      \put(7026,2290){\makebox(0,0)[l]{\strut{}$\theta$ = 1.1}}%
    }%
    \gplbacktext
    \put(0,0){\includegraphics{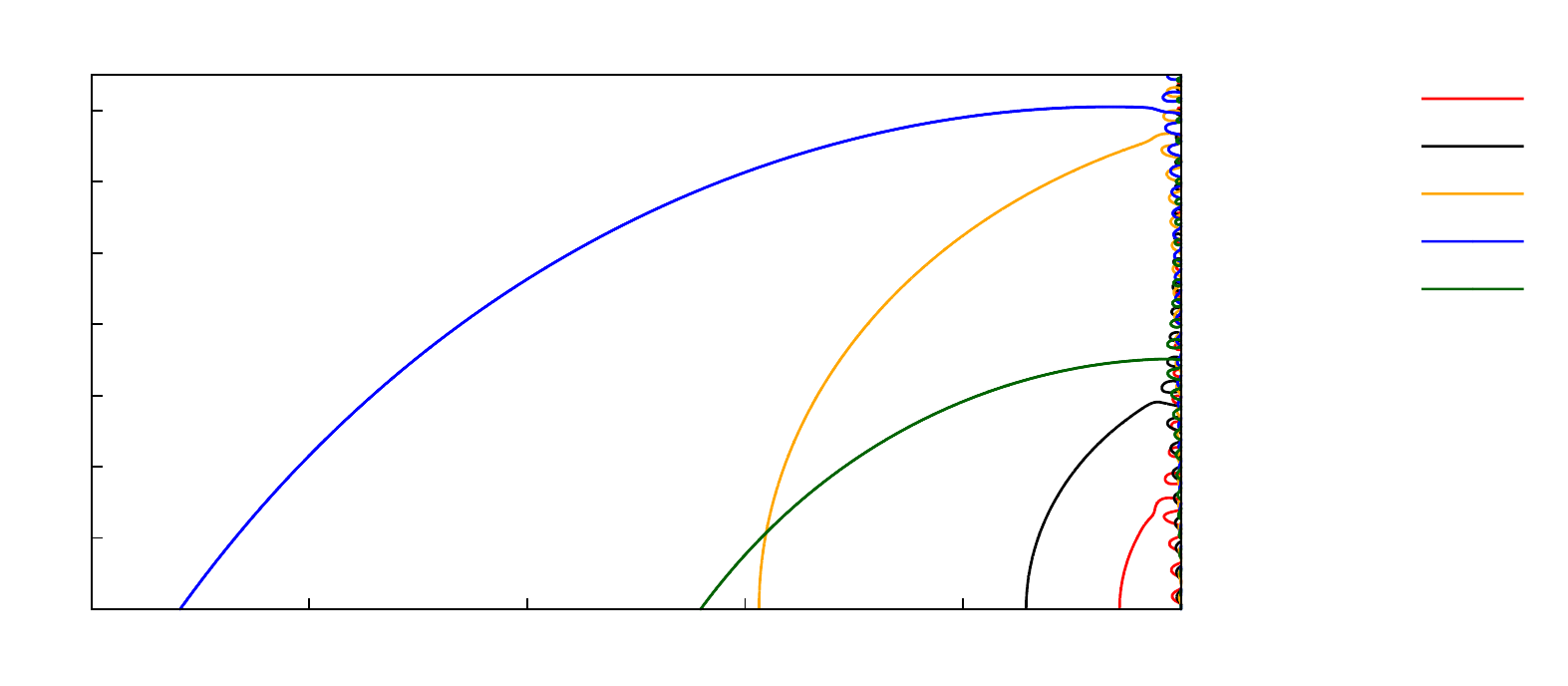}}%
    \gplfronttext
  \end{picture}%
  \endgroup
  \caption{Stability domain $D_\theta$ for the $\RL_3$ scheme. In the particular
    case $\theta=0$ (no stabilization, corresponding to the Adams-Bashforth
    scheme of order 3), the stability domain crosses the $x$-axis at
    $x\simeq-0.54$ (dark blue arrow).}
  \label{fig:RL3-D-stab}
\end{figure}

\begin{figure}
  \centering
  % GNUPLOT: LaTeX picture with Postscript
  \begingroup
  \sizeg
  \makeatletter
  \providecommand\color[2][]{%
    \GenericError{(gnuplot) \space\space\space\@spaces}{%
      Package color not loaded in conjunction with
      terminal option `colourtext'%
    }{See the gnuplot documentation for explanation.%
    }{Either use 'blacktext' in gnuplot or load the package
      color.sty in LaTeX.}%
    \renewcommand\color[2][]{}%
  }%
  \providecommand\includegraphics[2][]{%
    \GenericError{(gnuplot) \space\space\space\@spaces}{%
      Package graphicx or graphics not loaded%
    }{See the gnuplot documentation for explanation.%
    }{The gnuplot epslatex terminal needs graphicx.sty or graphics.sty.}%
    \renewcommand\includegraphics[2][]{}%
  }%
  \providecommand\rotatebox[2]{#2}%
  \@ifundefined{ifGPcolor}{%
    \newif\ifGPcolor
    \GPcolortrue
  }{}%
  \@ifundefined{ifGPblacktext}{%
    \newif\ifGPblacktext
    \GPblacktexttrue
  }{}%
  % define a \g@addto@macro without @ in the name:
  \let\gplgaddtomacro\g@addto@macro
  % define empty templates for all commands taking text:
  \gdef\gplbacktext{}%
  \gdef\gplfronttext{}%
  \makeatother
  \ifGPblacktext
  % no textcolor at all
  \def\colorrgb#1{}%
  \def\colorgray#1{}%
  \else
  % gray or color?
  \ifGPcolor
  \def\colorrgb#1{\color[rgb]{#1}}%
  \def\colorgray#1{\color[gray]{#1}}%
  \expandafter\def\csname LTw\endcsname{\color{white}}%
  \expandafter\def\csname LTb\endcsname{\color{black}}%
  \expandafter\def\csname LTa\endcsname{\color{black}}%
  \expandafter\def\csname LT0\endcsname{\color[rgb]{1,0,0}}%
  \expandafter\def\csname LT1\endcsname{\color[rgb]{0,1,0}}%
  \expandafter\def\csname LT2\endcsname{\color[rgb]{0,0,1}}%
  \expandafter\def\csname LT3\endcsname{\color[rgb]{1,0,1}}%
  \expandafter\def\csname LT4\endcsname{\color[rgb]{0,1,1}}%
  \expandafter\def\csname LT5\endcsname{\color[rgb]{1,1,0}}%
  \expandafter\def\csname LT6\endcsname{\color[rgb]{0,0,0}}%
  \expandafter\def\csname LT7\endcsname{\color[rgb]{1,0.3,0}}%
  \expandafter\def\csname LT8\endcsname{\color[rgb]{0.5,0.5,0.5}}%
  \else
  % gray
  \def\colorrgb#1{\color{black}}%
  \def\colorgray#1{\color[gray]{#1}}%
  \expandafter\def\csname LTw\endcsname{\color{white}}%
  \expandafter\def\csname LTb\endcsname{\color{black}}%
  \expandafter\def\csname LTa\endcsname{\color{black}}%
  \expandafter\def\csname LT0\endcsname{\color{black}}%
  \expandafter\def\csname LT1\endcsname{\color{black}}%
  \expandafter\def\csname LT2\endcsname{\color{black}}%
  \expandafter\def\csname LT3\endcsname{\color{black}}%
  \expandafter\def\csname LT4\endcsname{\color{black}}%
  \expandafter\def\csname LT5\endcsname{\color{black}}%
  \expandafter\def\csname LT6\endcsname{\color{black}}%
  \expandafter\def\csname LT7\endcsname{\color{black}}%
  \expandafter\def\csname LT8\endcsname{\color{black}}%
  \fi
  \fi
  \setlength{\unitlength}{0.0500bp}%
  \begin{picture}(9070.00,3968.00)%
    \gplgaddtomacro\gplbacktext{%
      \csname LTb\endcsname%
      \put(396,440){\makebox(0,0)[r]{\strut{} 0}}%
      \put(396,1057){\makebox(0,0)[r]{\strut{} 10}}%
      \put(396,1675){\makebox(0,0)[r]{\strut{} 20}}%
      \put(396,2292){\makebox(0,0)[r]{\strut{} 30}}%
      \put(396,2910){\makebox(0,0)[r]{\strut{} 40}}%
      \put(396,3527){\makebox(0,0)[r]{\strut{} 50}}%
      \put(731,220){\makebox(0,0){\strut{}-90}}%
      \put(1408,220){\makebox(0,0){\strut{}-80}}%
      \put(2085,220){\makebox(0,0){\strut{}-70}}%
      \put(2762,220){\makebox(0,0){\strut{}-60}}%
      \put(3440,220){\makebox(0,0){\strut{}-50}}%
      \put(4117,220){\makebox(0,0){\strut{}-40}}%
      \put(4794,220){\makebox(0,0){\strut{}-30}}%
      \put(5471,220){\makebox(0,0){\strut{}-20}}%
      \put(6148,220){\makebox(0,0){\strut{}-10}}%
      \put(6825,220){\makebox(0,0){\strut{} 0}}%
      \put(6805,440){\rotatebox{260}{\makebox(0,0)[l]{\strut{}\textcolor{NavyBlue} {$\leftarrow$ -0.3 }}}}%
      \put(6669,440){\rotatebox{240}{\makebox(0,0)[l]{\strut{}\textcolor{red} {$\leftarrow$ -2.3 }}}}%
      \put(6453,440){\rotatebox{100}{\makebox(0,0)[l]{\strut{}\textcolor{black}{$\leftarrow$ -5.5}}}}%
      \put(5065,440){\rotatebox{100}{\makebox(0,0)[l]{\strut{}\textcolor{Orange}{$\leftarrow$ -26.0}}}}%
      \put(968,440){\rotatebox{100}{\makebox(0,0)[l]{\strut{}\textcolor{blue}{$\leftarrow$ -86.5}}}}%
      \put(4110,440){\rotatebox{100}{\makebox(0,0)[l]{\strut{}\textcolor{OliveGreen}{$\leftarrow$ -40.1}}}}%
    }%
    \gplgaddtomacro\gplfronttext{%
      \csname LTb\endcsname%
      \put(7026,3390){\makebox(0,0)[l]{\strut{}$\theta$ = 0.85}}%
      \csname LTb\endcsname%
      \put(7026,3115){\makebox(0,0)[l]{\strut{}$\theta$ = 0.9}}%
      \csname LTb\endcsname%
      \put(7026,2840){\makebox(0,0)[l]{\strut{}$\theta$ = 0.95}}%
      \csname LTb\endcsname%
      \put(7026,2565){\makebox(0,0)[l]{\strut{}$\theta$ = 1.05}}%
      \csname LTb\endcsname%
      \put(7026,2290){\makebox(0,0)[l]{\strut{}$\theta$ = 1.1}}%
    }%
    \gplbacktext
    \put(0,0){\includegraphics{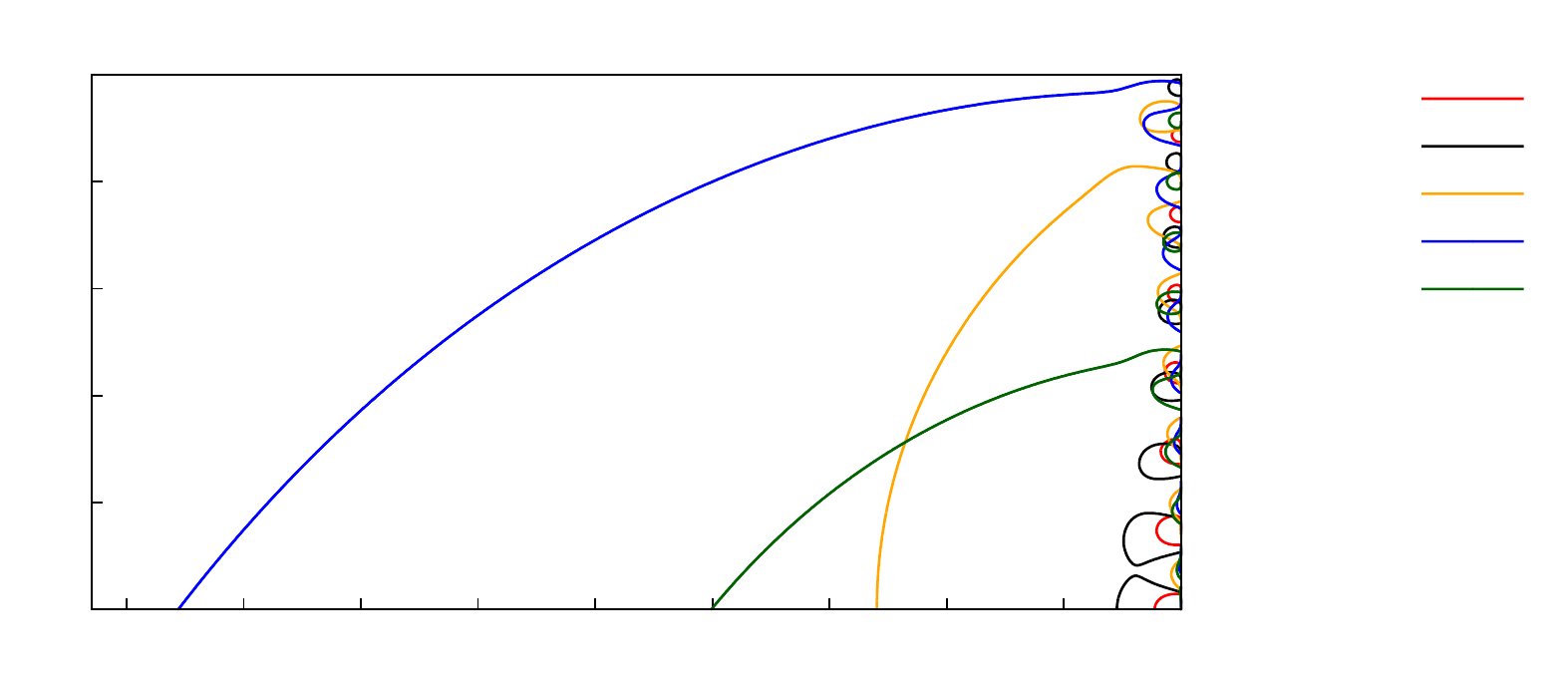}}%
    \gplfronttext
  \end{picture}%
  \endgroup
  \caption{Stability domain $D_\theta$ for the $\RL_4$ scheme. In the particular
    case $\theta=0$ (no stabilization, corresponding to the Adams-Bashforth
    scheme of order 3), the stability domain crosses the $x$-axis at
    $x\simeq-0.3$ (dark blue arrow).}
  \label{fig:RL4-D-stab}
\end{figure}

The situation is different for the Rush-Larsen methods of orders 3 and
4. The stability domains $D_\theta$ are
depicted on Figures~\ref{fig:RL3-D-stab} and \ref{fig:RL4-D-stab}, for various
values of $\theta$, and for the orders 3 and 4, respectively. Excepted
for the case $\theta=1$, the stability domain is always bounded: the scheme is
not $A(0)$-stable. However, the stability domain for $\theta\simeq 1$ is much
larger than without stabilization (corresponding to the Adams-Bashforth
schemes of orders 3 or 4). For the RL$_3$ scheme, the stability domain
when $\theta=0.85$ is 25 times wider on the left than $D_{\theta| \theta=0}$,
and when $\theta=1.05$ it is 400 times wider. For the RL$_4$ case,
$D_{\theta| \theta=1.05}$ is almost 300 times wider on the left than
$D_{\theta| \theta=0}$.
\section{Numerical results}
\label{sec:num-res}
In this section we present numerical experiments that illustrate the performances of the $\RL_k$ methods.
They will be compared to the exponential integrators of Adams type of order $k$ defined by equation (\ref{eab-scheme}), 
shortly denoted by $\EAB_k$.
The $\EAB_k$ schemes have been numerically studied 
in
\cite{coudiere-lontsi-pierre-2015},
for the resolution of the membrane equation in electrophysiology,
as compared to several classical methods.
In that context, they have been shown to be as stable as implicit methods with a much smaller cost.
We present the same numerical tests here,
so as to extend the comparison to the schemes benchmarked in \cite{coudiere-lontsi-pierre-2015}.
% In addition, several efficiency comparison in terms of CPU times are available in 
% \cite{Cari_2016}, for the membrane equation (\ref{eq:edo-ionic-w}) and
% for physiologically relevant quantities.

%% 
%% 
%% 
%% 
\subsection{The membrane equation}
\label{sec:memb-eq}
\begin{figure}[h!]
  \centering
  % GNUPLOT: LaTeX picture with Postscript
  \begingroup
  \sizeg
  \makeatletter
  \providecommand\color[2][]{%
    \GenericError{(gnuplot) \space\space\space\@spaces}{%
      Package color not loaded in conjunction with
      terminal option `colourtext'%
    }{See the gnuplot documentation for explanation.%
    }{Either use 'blacktext' in gnuplot or load the package
      color.sty in LaTeX.}%
    \renewcommand\color[2][]{}%
  }%
  \providecommand\includegraphics[2][]{%
    \GenericError{(gnuplot) \space\space\space\@spaces}{%
      Package graphicx or graphics not loaded%
    }{See the gnuplot documentation for explanation.%
    }{The gnuplot epslatex terminal needs graphicx.sty or graphics.sty.}%
    \renewcommand\includegraphics[2][]{}%
  }%
  \providecommand\rotatebox[2]{#2}%
  \@ifundefined{ifGPcolor}{%
    \newif\ifGPcolor
    \GPcolorfalse
  }{}%
  \@ifundefined{ifGPblacktext}{%
    \newif\ifGPblacktext
    \GPblacktexttrue
  }{}%
  % define a \g@addto@macro without @ in the name:
  \let\gplgaddtomacro\g@addto@macro
  % define empty templates for all commands taking text:
  \gdef\gplbacktext{}%
  \gdef\gplfronttext{}%
  \makeatother
  \ifGPblacktext
  % no textcolor at all
  \def\colorrgb#1{}%
  \def\colorgray#1{}%
  \else
  % gray or color?
  \ifGPcolor
  \def\colorrgb#1{\color[rgb]{#1}}%
  \def\colorgray#1{\color[gray]{#1}}%
  \expandafter\def\csname LTw\endcsname{\color{white}}%
  \expandafter\def\csname LTb\endcsname{\color{black}}%
  \expandafter\def\csname LTa\endcsname{\color{black}}%
  \expandafter\def\csname LT0\endcsname{\color[rgb]{1,0,0}}%
  \expandafter\def\csname LT1\endcsname{\color[rgb]{0,1,0}}%
  \expandafter\def\csname LT2\endcsname{\color[rgb]{0,0,1}}%
  \expandafter\def\csname LT3\endcsname{\color[rgb]{1,0,1}}%
  \expandafter\def\csname LT4\endcsname{\color[rgb]{0,1,1}}%
  \expandafter\def\csname LT5\endcsname{\color[rgb]{1,1,0}}%
  \expandafter\def\csname LT6\endcsname{\color[rgb]{0,0,0}}%
  \expandafter\def\csname LT7\endcsname{\color[rgb]{1,0.3,0}}%
  \expandafter\def\csname LT8\endcsname{\color[rgb]{0.5,0.5,0.5}}%
  \else
  % gray
  \def\colorrgb#1{\color{black}}%
  \def\colorgray#1{\color[gray]{#1}}%
  \expandafter\def\csname LTw\endcsname{\color{white}}%
  \expandafter\def\csname LTb\endcsname{\color{black}}%
  \expandafter\def\csname LTa\endcsname{\color{black}}%
  \expandafter\def\csname LT0\endcsname{\color{black}}%
  \expandafter\def\csname LT1\endcsname{\color{black}}%
  \expandafter\def\csname LT2\endcsname{\color{black}}%
  \expandafter\def\csname LT3\endcsname{\color{black}}%
  \expandafter\def\csname LT4\endcsname{\color{black}}%
  \expandafter\def\csname LT5\endcsname{\color{black}}%
  \expandafter\def\csname LT6\endcsname{\color{black}}%
  \expandafter\def\csname LT7\endcsname{\color{black}}%
  \expandafter\def\csname LT8\endcsname{\color{black}}%
  \fi
  \fi
  \setlength{\unitlength}{0.0500bp}%
  \begin{picture}(3684.00,2550.00)%
    \gplgaddtomacro\gplbacktext{%
      \csname LTb\endcsname%
      \put(264,440){\makebox(0,0)[r]{\strut{}-100}}%
      \put(264,741){\makebox(0,0)[r]{\strut{}-80}}%
      \put(264,1043){\makebox(0,0)[r]{\strut{}-60}}%
      \put(264,1344){\makebox(0,0)[r]{\strut{}-40}}%
      \put(264,1645){\makebox(0,0)[r]{\strut{}-20}}%
      \put(264,1946){\makebox(0,0)[r]{\strut{} 0}}%
      \put(264,2248){\makebox(0,0)[r]{\strut{} 20}}%
      \put(264,2549){\makebox(0,0)[r]{\strut{} 40}}%
      \put(396,220){\makebox(0,0){\strut{} 0}}%
      \put(1152,220){\makebox(0,0){\strut{} 100}}%
      \put(1908,220){\makebox(0,0){\strut{} 200}}%
      \put(2663,220){\makebox(0,0){\strut{} 300}}%
      \put(3419,220){\makebox(0,0){\strut{} 400}}%
      \put(-413,1494){\rotatebox{-270}{\makebox(0,0){\strut{}Potential (\si{mV})}}}%
      \put(1907,-110){\makebox(0,0){\strut{}Time (\si{ms})}}%
      \put(1907,2879){\makebox(0,0){\strut{}Transmembrane potential}}%
    }%
    \gplgaddtomacro\gplfronttext{%
      \csname LTb\endcsname%
      \put(2432,2376){\makebox(0,0)[r]{\strut{}$v$}}%
    }%
    \gplbacktext
    \put(0,0){\includegraphics{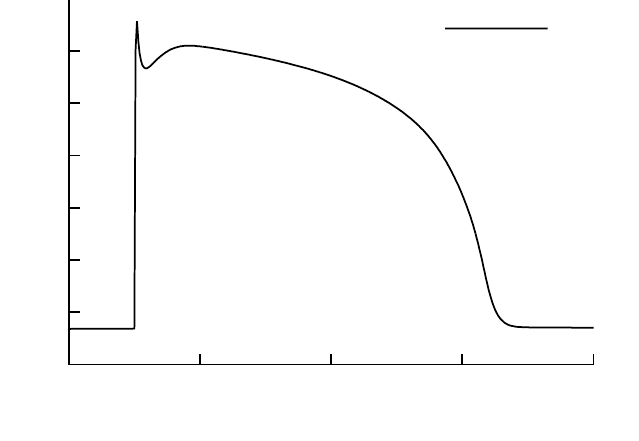}}%
    \gplfronttext
  \end{picture}%
  \endgroup
  $\quad \quad \quad \quad $
  % GNUPLOT: LaTeX picture with Postscript
  \begingroup
  \sizeg
  \makeatletter
  \providecommand\color[2][]{%
    \GenericError{(gnuplot) \space\space\space\@spaces}{%
      Package color not loaded in conjunction with
      terminal option `colourtext'%
    }{See the gnuplot documentation for explanation.%
    }{Either use 'blacktext' in gnuplot or load the package
      color.sty in LaTeX.}%
    \renewcommand\color[2][]{}%
  }%
  \providecommand\includegraphics[2][]{%
    \GenericError{(gnuplot) \space\space\space\@spaces}{%
      Package graphicx or graphics not loaded%
    }{See the gnuplot documentation for explanation.%
    }{The gnuplot epslatex terminal needs graphicx.sty or graphics.sty.}%
    \renewcommand\includegraphics[2][]{}%
  }%
  \providecommand\rotatebox[2]{#2}%
  \@ifundefined{ifGPcolor}{%
    \newif\ifGPcolor
    \GPcolorfalse
  }{}%
  \@ifundefined{ifGPblacktext}{%
    \newif\ifGPblacktext
    \GPblacktexttrue
  }{}%
  % define a \g@addto@macro without @ in the name:
  \let\gplgaddtomacro\g@addto@macro
  % define empty templates for all commands taking text:
  \gdef\gplbacktext{}%
  \gdef\gplfronttext{}%
  \makeatother
  \ifGPblacktext
  % no textcolor at all
  \def\colorrgb#1{}%
  \def\colorgray#1{}%
  \else
  % gray or color?
  \ifGPcolor
  \def\colorrgb#1{\color[rgb]{#1}}%
  \def\colorgray#1{\color[gray]{#1}}%
  \expandafter\def\csname LTw\endcsname{\color{white}}%
  \expandafter\def\csname LTb\endcsname{\color{black}}%
  \expandafter\def\csname LTa\endcsname{\color{black}}%
  \expandafter\def\csname LT0\endcsname{\color[rgb]{1,0,0}}%
  \expandafter\def\csname LT1\endcsname{\color[rgb]{0,1,0}}%
  \expandafter\def\csname LT2\endcsname{\color[rgb]{0,0,1}}%
  \expandafter\def\csname LT3\endcsname{\color[rgb]{1,0,1}}%
  \expandafter\def\csname LT4\endcsname{\color[rgb]{0,1,1}}%
  \expandafter\def\csname LT5\endcsname{\color[rgb]{1,1,0}}%
  \expandafter\def\csname LT6\endcsname{\color[rgb]{0,0,0}}%
  \expandafter\def\csname LT7\endcsname{\color[rgb]{1,0.3,0}}%
  \expandafter\def\csname LT8\endcsname{\color[rgb]{0.5,0.5,0.5}}%
  \else
  % gray
  \def\colorrgb#1{\color{black}}%
  \def\colorgray#1{\color[gray]{#1}}%
  \expandafter\def\csname LTw\endcsname{\color{white}}%
  \expandafter\def\csname LTb\endcsname{\color{black}}%
  \expandafter\def\csname LTa\endcsname{\color{black}}%
  \expandafter\def\csname LT0\endcsname{\color{black}}%
  \expandafter\def\csname LT1\endcsname{\color{black}}%
  \expandafter\def\csname LT2\endcsname{\color{black}}%
  \expandafter\def\csname LT3\endcsname{\color{black}}%
  \expandafter\def\csname LT4\endcsname{\color{black}}%
  \expandafter\def\csname LT5\endcsname{\color{black}}%
  \expandafter\def\csname LT6\endcsname{\color{black}}%
  \expandafter\def\csname LT7\endcsname{\color{black}}%
  \expandafter\def\csname LT8\endcsname{\color{black}}%
  \fi
  \fi
  \setlength{\unitlength}{0.0500bp}%
  \begin{picture}(3684.00,2550.00)%
    \gplgaddtomacro\gplbacktext{%
      \csname LTb\endcsname%
      \put(396,440){\makebox(0,0)[r]{\strut{}-200}}%
      \put(396,1319){\makebox(0,0)[r]{\strut{}-100}}%
      \put(396,2198){\makebox(0,0)[r]{\strut{} 0}}%
      \put(528,220){\makebox(0,0){\strut{} 0}}%
      \put(1218,220){\makebox(0,0){\strut{} 100}}%
      \put(1908,220){\makebox(0,0){\strut{} 200}}%
      \put(2597,220){\makebox(0,0){\strut{} 300}}%
      \put(3287,220){\makebox(0,0){\strut{} 400}}%
      \put(-281,1494){\rotatebox{-270}{\makebox(0,0){\strut{}Current (\si{A/F/cm^2})}}}%
      \put(1907,-110){\makebox(0,0){\strut{}Time (\si{ms})}}%
      \put(1907,2879){\makebox(0,0){\strut{}Fast inward sodium current}}%
    }%
    \gplgaddtomacro\gplfronttext{%
      \csname LTb\endcsname%
      \put(2300,2376){\makebox(0,0)[r]{\strut{}$I_{Na}$}}%
    }%
    \gplbacktext
    \put(0,0){\includegraphics{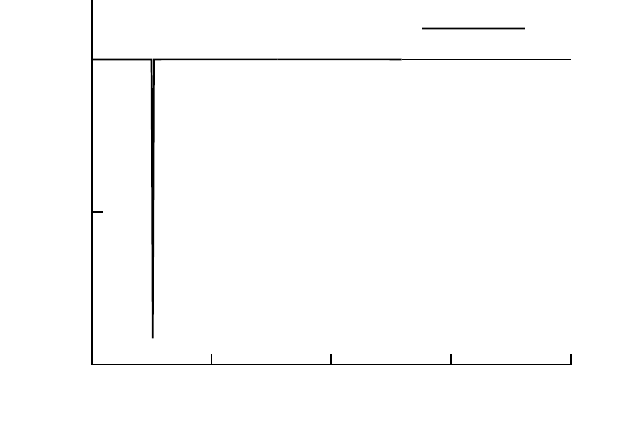}}%
    \gplfronttext
  \end{picture}%
  \endgroup
  \caption{TNNP model illustration. Left, cellular 
    action potential: starting from a (negative) rest value, the transmembrane potential $v(t)$ has a stiff depolarization followed by a plateau and a repolarization to the rest value. Right, depolarization is induced by an ionic sodium current $I_{Na}$, with obvious large stiffness.% of which is illustrated here at action potential duration time scale.
  }
  \label{fig:BR}
\end{figure}
% 
% We consider a class of models in cardiac electrophysiology.
The cellular action potential for cardiac cells is described on Figure \ref{fig:BR}.
This phenomenon  displays a stiff behavior characterized by the presence of heterogeneous time scales.
The electrical activity of cardiac cells is modeled
with an ODE system called membrane equation. It has the form 
\begin{align}
  \label{eq:edo-ionic-w}
  \diff{v}{t}  = - I_{ion}(v,w,c) + I_{st}(t),\quad 
  \diff{w_i}{t} = \frac{w_{\infty,i}(v) - w_i}{\tau_i(v)} 
  ,\qquad 
  \diff{c}{t} = g(v,w,c), 
\end{align}
where $w=(w_1,\dots,w_p)\in\R^p$ is a vector of  gating
variables, $c\in \R^q$ is a vector of ionic concentrations, and $v\in\R$ is the transmembrane potential, 
we refer to \cite{Hodgkin52,beeler-reuter,LR2a,tnnp} for details.
The gating variables describe the opening state (between 0 and 1) of various protein structures on the cell membrane, which control ionic transfers between the intra and extra-cellular media.
Each gating variable $w_i$
evolves towards the state $w_{\infty,i}(v)$ at rate $\tau_i(v)$.
Specific ionic currents (sodium, potassium, ...) across the cellular membrane are computed with the help of the variables $v$, $w$ and $c$. 
The sum of these currents
% With the help of these variables are modeled ionic currents for specific ionic species (sodium, potassium, ...) the sum of which 
defines the total ionic current  
$I_{ion}(v,w,c)$ across the membrane.
The function
$I_{st}(t)$ is a  source term, it represents a stimulation current. 
The membrane equation corresponds to the ODE system in the monodomain model (\ref{F0}) with $\zeta=(w,c)$.

We will consider two such models: the BR model \cite{beeler-reuter}
and the TNNP model \cite{tnnp}.% for human cardiac cells.
The BR model \cite{beeler-reuter}
describes the membrane action potential of mammalian ventricular myocardial cells.
It involves 6 gating variables ($p=6$, 
they are denoted by $m,~h~,j~,d~,f~,x_i$) and one ionic concentration ($q=1$): the intra-cellular calcium $[Ca_i]$.
% Each gating variables is controlled by an ODE as in equation \eqref{eq:edo-ionic-w}.
The Nernst potential $E_{Ca}$ for the calcium ions then is time-dependent and a (slow inward) calcium current $I_s$ is modeled as
$I_s = g_s d f (v-E_{Ca})$ depending on the gating variables $d$, $f$, the transmembrane potential $v$ and a constant $g_s$.
A fast inward sodium current $I_{Na}$ (depicted on Figure \ref{fig:BR}) that depends on the three gating variables $m$, $h$, $j$ and on $v$ is similarly described.
Two outward currents are modeled: $I_{x_1}$ that depends on $x_1$ and $v$ and $I_{K_1}$ that only depends on $v$.
The total ionic current in (\ref{eq:edo-ionic-w}) is the sum of these four currents
$I_{Na} + I_{s} + I_{K_1} + I_{x_1} = -I_{ion}(w,c,v)$. 
\\
The TNNP model is specifically designed for human ventricular myocytes.
It is more sophisticated than the BR model but conserves the same general structure.
It involves 12 gating variables and 4 ionic concentrations ($p=12$ and $q=4$). The total ionic current $I_{ion}$ is the sum of 15 specific ionic currents.
\\ \\
The membrane equation (\ref{eq:edo-ionic-w}) can be reformulated in the form of  (\ref{F2}) with
\begin{equation}
  \label{eq:edo-ionic}
  a(t,y) = 
  \begin{pmatrix}
    -1/\tau(v) & 0& 0 \\
    0 & 0& 0 \\
    0 & 0& 0
  \end{pmatrix},
  \quad b(t,y) =
  \begin{pmatrix}
    w_\infty(v)/\tau(v) \\
    g(y) \\
    -I_{ion}(y) + I_{st}(t)
  \end{pmatrix},
\end{equation}
for $y = (w,c,v) \in \R^N$ ($N=p+q+1$) and where $-1/\tau(v)$ is the
$p \times p$ diagonal matrix with diagonal entries
$\left(-1/\tau_i(v)\right)_{i=1,\ldots, p}$.
The resulting matrix $a(t,y)$ is diagonal.

\subsection{Convergence}  

No analytical solution is 
available for the chosen application.
A reference solution $y_{\text{ref}}$ for a reference  time step $h_{\text{ref}}$ is computed with the Runge-Kutta scheme of order 4 to analyze the convergence properties of the $\RL_k$ schemes.
Numerical solutions $y$ are compared to $y_{\text{ref}}$ for coarsest time steps $h = 2^m h_{\text{ref}}$.% for increasing $m$.
\\
A numerical solution $y$ consists in successive values $y_n$ at the time instants $t_n = n h$.
On every interval $(t_{3n}, t_{3n+3})$ the polynomial $\overline{y}$  of degree at most 3 so that 
$\overline{y}(t_{3n+i}) = y_{3n+i}$, 
$i=0,\ldots, 3$
is constructed.
On $(0,T)$ $\overline{y}$ is continuous and piecewise  polynomial of degree 3\, its values at the reference time instants $n h_{\text{ref}}$ are computed. This provides a projection $P(y)$ of the numerical solution $y$ onto the reference grid.
Then $P(y)$ can be compared with the reference solution $y_{\text{ref}}$.
The numerical error is defined by
\begin{equation}
  \label{eq:def-erreur}
  e(h)= \frac{ \max \left|v_{\text{ref}} - P(v)\right|}{ \max \left|v_{\text{ref}}\right|},
\end{equation}
where the potential $v$ is the last and stiffest component of $y$ in equation \eqref{eq:edo-ionic-w}.
\begin{figure}[h!]
  % GNUPLOT: LaTeX picture with Postscript
  \begingroup
  \sizeg
  \makeatletter
  \providecommand\color[2][]{%
    \GenericError{(gnuplot) \space\space\space\@spaces}{%
      Package color not loaded in conjunction with
      terminal option `colourtext'%
    }{See the gnuplot documentation for explanation.%
    }{Either use 'blacktext' in gnuplot or load the package
      color.sty in LaTeX.}%
    \renewcommand\color[2][]{}%
  }%
  \providecommand\includegraphics[2][]{%
    \GenericError{(gnuplot) \space\space\space\@spaces}{%
      Package graphicx or graphics not loaded%
    }{See the gnuplot documentation for explanation.%
    }{The gnuplot epslatex terminal needs graphicx.sty or graphics.sty.}%
    \renewcommand\includegraphics[2][]{}%
  }%
  \providecommand\rotatebox[2]{#2}%
  \@ifundefined{ifGPcolor}{%
    \newif\ifGPcolor
    \GPcolorfalse
  }{}%
  \@ifundefined{ifGPblacktext}{%
    \newif\ifGPblacktext
    \GPblacktexttrue
  }{}%
  % define a \g@addto@macro without @ in the name:
  \let\gplgaddtomacro\g@addto@macro
  % define empty templates for all commands taking text:
  \gdef\gplbacktext{}%
  \gdef\gplfronttext{}%
  \makeatother
  \ifGPblacktext
  % no textcolor at all
  \def\colorrgb#1{}%
  \def\colorgray#1{}%
  \else
  % gray or color?
  \ifGPcolor
  \def\colorrgb#1{\color[rgb]{#1}}%
  \def\colorgray#1{\color[gray]{#1}}%
  \expandafter\def\csname LTw\endcsname{\color{white}}%
  \expandafter\def\csname LTb\endcsname{\color{black}}%
  \expandafter\def\csname LTa\endcsname{\color{black}}%
  \expandafter\def\csname LT0\endcsname{\color[rgb]{1,0,0}}%
  \expandafter\def\csname LT1\endcsname{\color[rgb]{0,1,0}}%
  \expandafter\def\csname LT2\endcsname{\color[rgb]{0,0,1}}%
  \expandafter\def\csname LT3\endcsname{\color[rgb]{1,0,1}}%
  \expandafter\def\csname LT4\endcsname{\color[rgb]{0,1,1}}%
  \expandafter\def\csname LT5\endcsname{\color[rgb]{1,1,0}}%
  \expandafter\def\csname LT6\endcsname{\color[rgb]{0,0,0}}%
  \expandafter\def\csname LT7\endcsname{\color[rgb]{1,0.3,0}}%
  \expandafter\def\csname LT8\endcsname{\color[rgb]{0.5,0.5,0.5}}%
  \else
  % gray
  \def\colorrgb#1{\color{black}}%
  \def\colorgray#1{\color[gray]{#1}}%
  \expandafter\def\csname LTw\endcsname{\color{white}}%
  \expandafter\def\csname LTb\endcsname{\color{black}}%
  \expandafter\def\csname LTa\endcsname{\color{black}}%
  \expandafter\def\csname LT0\endcsname{\color{black}}%
  \expandafter\def\csname LT1\endcsname{\color{black}}%
  \expandafter\def\csname LT2\endcsname{\color{black}}%
  \expandafter\def\csname LT3\endcsname{\color{black}}%
  \expandafter\def\csname LT4\endcsname{\color{black}}%
  \expandafter\def\csname LT5\endcsname{\color{black}}%
  \expandafter\def\csname LT6\endcsname{\color{black}}%
  \expandafter\def\csname LT7\endcsname{\color{black}}%
  \expandafter\def\csname LT8\endcsname{\color{black}}%
  \fi
  \fi
  \setlength{\unitlength}{0.0500bp}%
  \begin{picture}(7936.00,3400.00)%
    \gplgaddtomacro\gplbacktext{%
      \csname LTb\endcsname%
      \put(660,669){\makebox(0,0)[r]{\strut{} 1e-10}}%
      \csname LTb\endcsname%
      \put(660,1127){\makebox(0,0)[r]{\strut{} 1e-08}}%
      \csname LTb\endcsname%
      \put(660,1585){\makebox(0,0)[r]{\strut{} 1e-06}}%
      \csname LTb\endcsname%
      \put(660,2043){\makebox(0,0)[r]{\strut{} 0.0001}}%
      \csname LTb\endcsname%
      \put(660,2501){\makebox(0,0)[r]{\strut{} 0.01}}%
      \csname LTb\endcsname%
      \put(660,2959){\makebox(0,0)[r]{\strut{} 1}}%
      \csname LTb\endcsname%
      \put(1801,220){\makebox(0,0){\strut{} 0.001}}%
      \csname LTb\endcsname%
      \put(3245,220){\makebox(0,0){\strut{} 0.01}}%
      \csname LTb\endcsname%
      \put(4690,220){\makebox(0,0){\strut{} 0.1}}%
      \put(-110,1699){\makebox(0,0){\strut{} $e(h)$}}%
      \put(3175,-110){\makebox(0,0){\strut{}Time step}}%
    }%
    \gplgaddtomacro\gplfronttext{%
      \csname LTb\endcsname%
      \put(5760,2822){\makebox(0,0)[l]{\strut{}$RL2$}}%
      \csname LTb\endcsname%
      \put(5760,2547){\makebox(0,0)[l]{\strut{}$RL3$}}%
      \csname LTb\endcsname%
      \put(5760,2272){\makebox(0,0)[l]{\strut{}$RL4$}}%
      \csname LTb\endcsname%
      \put(5760,1997){\makebox(0,0)[l]{\strut{} slope $2 $}}%
      \csname LTb\endcsname%
      \put(5760,1722){\makebox(0,0)[l]{\strut{} slope $3 $}}%
      \csname LTb\endcsname%
      \put(5760,1447){\makebox(0,0)[l]{\strut{} slope $4 $}}%
    }%
    \gplbacktext
    \put(0,0){\includegraphics{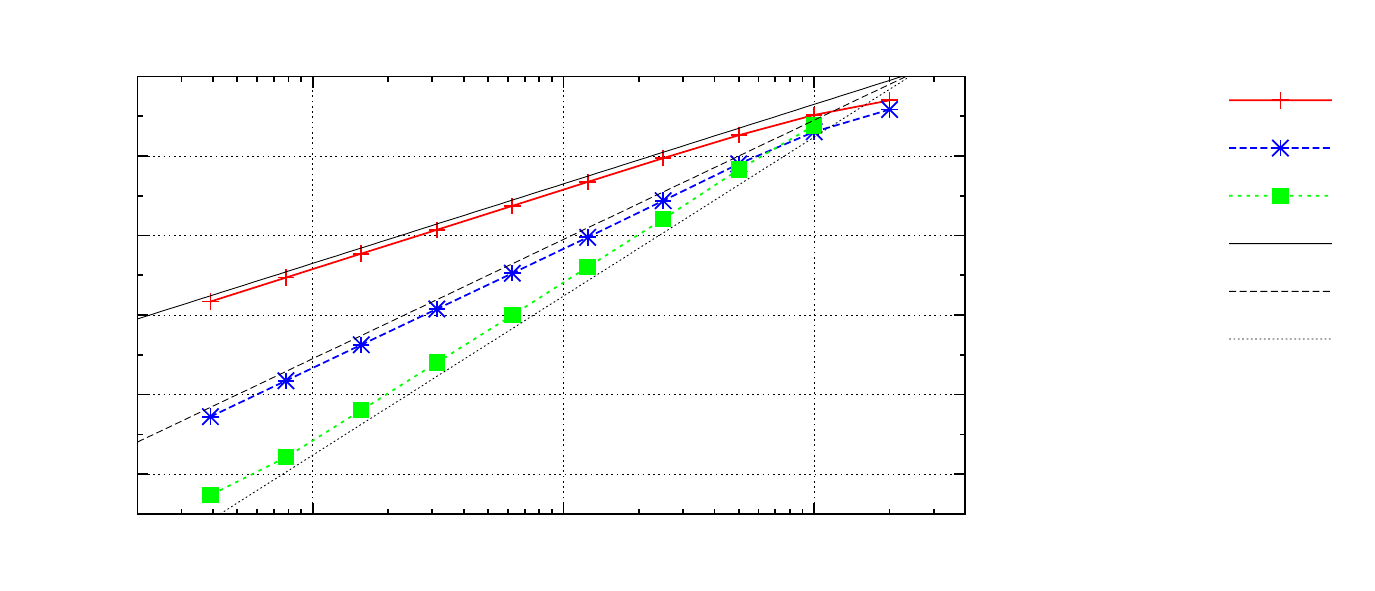}}%
    \gplfronttext
  \end{picture}%
  \endgroup
  \caption{
    % \blue{Figure : il faut mettre les trois pentes, et certainement adapter l'ÃÂ©chelle en y pour que le tout reste lisible. On peut aussi faire deux courbes cÃÂ´te ÃÂ  cÃÂ´te. }
    Relative error $e(h)$ (definition \eqref{eq:def-erreur}) as a function of the time step $h$ for the $\RL_k$ schemes, for $k=2$ to $4$ and in Log/Log scale.} 
  \label{fig:The numerical convergence}
\end{figure}
The convergence graphs for the BR model are plotted on Figure \ref{fig:The numerical convergence}. Each scheme displays the
expected asymptotic behavior of 
Proposition \ref{prop:stab}: $e(h)=O(h^k)$ as $h\to 0$.
\subsection{Stability}

Spiteri et al. in \cite{spiteri} have evaluated the stiffness of the BR and TNNP models along one cellular action potential (as depicted on Figure \ref{fig:BR}).
The largest negative real  part of the eigenvalues of the Jacobian matrix during the action potential is of $-1170$ and $-82$ for the TNNP and BR models, respectively. 
The TNNP model thus is 15 times stiffer than
the BR model ($15 \simeq 1170/82$).

Robustness to stiffness for the $\RL_k$ schemes is evaluated by comparing the critical time  steps for these two models. The critical time  step  $\Delta t_0$ is defined as the largest time step such that  the numerical simulations run 
without overflow for $h< \Delta t_0$.
The results are presented in Table \ref{tab_crit_tim1}.
\begin{table}[!h]
  \centering
  \caption{
    Critical time steps $\Delta t_0$ for the $\RL_k$ and $\EAB_k$ schemes}
  \label{tab_crit_tim1}
  \begin{tabular}{c|c|c|c||c|c|c}
    method & $\RL_2$ & $\RL_3$ & $\RL_4$& $\EAB_2$  & $\EAB_3$ & $\EAB_4$ \\ 
    \hline
    BR & 0.323 & 0.200&  0.149  & 0.424 & 0.203 & 0.123     \\ 
    \hline
    TNNP & 0.120 & 0.148 & 0.111& 0.233 & 0.108 &7.56 $10^{-2}$ \\ 
  \end{tabular}
\end{table}
\\
An excellent robustness to stiffness can be observed.
the $\RL_k$ schemes are not $A(\alpha)$ stable, and the
critical time step
is expected 
to be divided by 15 in case of an increase of stiffness of magnitude 15.
It is here divided by 2.7, 2.0 and 1.3 for $k=2$, 3 and 4, respectively.
A comparison with the $\EAB_k$ schemes shows that the two schemes 
have similar 
robustness to stiffness. 
% In conjunction with the numerical results in \cite{coudiere-lontsi-pierre-2015}, this means that the $\RL_k$ scheme robustness to stiffness is as efficient as for the $BDF_k$ implicit methods.
% This is a remarkable property for an explicit method.
% \\
% For a method that is not $A(\alpha)$ stable,  
% This is not observed here, though 
% \blue{The reason  is} that the ODE system (\ref{eq:edo-ionic-w}) is only partially stabilized by (\ref{eq:edo-ionic}).
Loss of stability is induced by the non-stabilized part, whose eigenvalues are less modified by the change of model.

\subsection{Accuracy}
The $\RL_k$ schemes are  compared to the $\EAB_k$ schemes in terms of accuracy.
This  is done using  the relative error $e(h)$ in equation (\ref{eq:def-erreur}), for the BR and TNNP models (we recall than the TNNP model is stiffer by a factor of 15). The results are collected in Tables \ref{acc_tim1} and \ref{acc_tim2}.
\begin{table}[!h]
  \centering
  \caption{Relative error $e(h)$ (eq. \eqref{eq:def-erreur}) for the BR model.}
  \label{acc_tim1}
  \begin{tabular}{c|c|c|c||c|c|c}
    $h$       & $ \RL_2$     & $\RL_3$     & $\RL_4$      & $\EAB_2$     & $\EAB_3$     & $\EAB_4$ \\ 
    \hline
    0.2       & 0.251 & 0.147 & -           & 0.284 & 0.516 & -            \\ 
    \hline
    0.1       & 0.107 & 4.07 $10^{-2}$& 5.86 $10^{-2}$& 9.26 $10^{-2}$& 9.17 $10^{-2}$& 0.119 \\ 
    \hline
    0.05 & 3.35 $10^{-2}$& 6.34 $10^{-3}$& 4.58 $10^{-3}$& 2.31 $10^{-2}$ & 1.09 $10^{-2}$& 8.96 $10^{-3}$\\ 
    \hline
    0.025 & 8.88 $10^{-3}$& 7.57 $10^{-4}$& 2.61 $10^{-4}$& 5.39 $10^{-3}$& 1.17 $10^{-3}$& 4.33 $10^{-4}$\\ 
  \end{tabular}
\end{table}
\begin{table}[!h]
  \centering
  \caption{Relative error $e(h)$ (eq. \eqref{eq:def-erreur}) for the TNNP model.}
  \label{acc_tim2}
  \begin{tabular}{c|c|c|c||c|c|c}
    $h$     & $ \RL_2$        & $\RL_3$          & $\RL_4$          & $\EAB_2$        & $\EAB_3$       & $\EAB_4$         \\ 
    \hline
    0.1     & 0.177         & 0.305          & 0.421          & 0.351          & 0.530          & -         \\ 
    \hline
    0.05    &7.39 $10^{-2}$ & 4.54 $10^{-2}$ & 4.61 $10^{-2}$ & 9.01 $10^{-2}$ & 5.59 $10^{-2}$ & 8.93 $10^{-2}$\\ 
    \hline
    0.025   & 2.21 $10^{-2}$& 6.53 $10^{-3}$ & 5.96 $10^{-3}$ & 2.14 $10^{-2}$ & 7.34 $10^{-3}$ & 8.34 $10^{-3}$\\ 
    \hline 
    0.0125  & 5.75 $10^{-3}$& 8.05 $10^{-4}$ &3.21 $10^{-4}$ & 5.11 $10^{-3}$ & 7.62 $10^{-4}$ & 3.70 $10^{-4}$\\ 
  \end{tabular}
\end{table}
\\
For the $\RL_2$ and the $\EAB_2$ schemes, the accuracies are very close, the $\EAB_2$ scheme being slightly more accurate for the BR model.
For the orders 3 and 4, 
% a \blue{non-negligible} difference is observed between the \blue{$\RL_k$ and the $\EAB_k$} schemes. 
the $\RL_k$ schemes are more accurate at large 
time steps. 
For smaller time steps, accuracies are almost the same.
The $\RL_k$ and $\EAB_k$ have the same accuracy in the asymptotic convergence region.
% Outside this region, \blue{the $\RL_k$ schemes are} more precise.

\section{Conclusion}
\label{sec:conclusion}
In this paper, we have introduced two new ODE solvers, that we have
called Rush-Larsen schemes of orders 3 and 4. They are explicit
multistep exponential integrators. Their definition is simple 
inducing an easy  implementation. 
We exposed the  analysis of convergence and of stability under perturbation  for these two schemes.
We also analyzed their Dahlquist stability: they are
not $A(0)$ stable, but exhibit very large stability domains for
sufficiently accurate stabilization. The numerical behavior of
the schemes is analyzed for a complex and realistic stiff application.
The $\RL_k$ schemes are as stable as exponential integrators of Adams type,
allowing simulations at large time step. 
On the presented example, 
the $\RL_k$ schemes are more accurate for $k=3$ and 4
than the exponential integrators of Adams type, 
when considering larger time steps.  
They are also shown to be robust to stiffness in terms of both stability and
accuracy.

\bibliographystyle{abbrv}
\bibliography{biblio}

\end{document}